\documentclass[nonblindrev]{informs3} 

\usepackage[english]{babel}
\usepackage{graphicx}
\usepackage{tabularx}
\usepackage{multirow}
\usepackage{hyperref}
\usepackage{array}
\usepackage{caption}
\usepackage{booktabs}
\usepackage{lscape}
\usepackage{longtable}
\usepackage{subcaption}
\usepackage{tablefootnote}
\usepackage{enumitem}

\usepackage{import}
\usepackage{float} 
\usepackage[flushleft]{threeparttable}
\usepackage{multicol}
\usepackage{amssymb}
\usepackage{tikz}

\usepackage[ruled, linesnumbered, noend]{algorithm2e}

\SetCommentSty{mycommfont}

\newcommand{\mli}[1]{\mathit{#1}}
\DoubleSpacedXI

\usepackage{natbib}
 \bibpunct[, ]{(}{)}{,}{a}{}{,}%
 \def\bibfont{\small}%
\TheoremsNumberedThrough     

\EquationsNumberedThrough    

\begin{document}


\RUNAUTHOR{Charaf et al.} 

\RUNTITLE{Branch-and-Price for the 2E-IRP}

\TITLE{A Branch-and-Price Algorithm for the Two-Echelon Inventory-Routing Problem}

\ARTICLEAUTHORS{%
\AUTHOR{Sara Charaf}
\AFF{School of Industrial Engineering, Eindhoven University of Technology, 5600MB Eindhoven, Netherlands, \EMAIL{s.charaf@tue.nl}, \URL{https://orcid.org/0000-0001-8100-8140}}
\AUTHOR{Duygu Ta\c s}
\AFF{Department of Industrial Engineering, MEF University, Istanbul, Turkey, \EMAIL{tasd@mef.edu.tr}, \URL{https://orcid.org/0000-0002-3579-4600}}
\AUTHOR{Simme Douwe P. Flapper, Tom Van Woensel}
\AFF{School of Industrial Engineering, Eindhoven University of Technology, 5600MB Eindhoven, Netherlands, \EMAIL{s.d.p.flapper@tue.nl}, \EMAIL{t.v.woensel@tue.nl}, \URL{https://orcid.org/
0000-0001-7559-3419}, \URL{https://orcid.org/0000-0003-4766-2346}}

} 

\ABSTRACT{
The two-echelon inventory-routing problem (2E-IRP) addresses the coordination of inventory management and freight transportation throughout a two-echelon supply network. 
The latter consists of geographically widespread customers whose demand over a discrete planning horizon can be met from either their local inventory or intermediate facilities' inventory. 
Intermediate facilities are located in the city outskirt and are supplied from distant suppliers. 
The 2E-IRP aims to minimize transportation costs and inventory costs while meeting customers' demand. 
A route-based formulation is proposed and a branch-and-price algorithm is developed for solving the 2E-IRP. 
A labeling algorithm is used to solve several pricing subproblems associated with each period and intermediate facility.
We generate 400 instances and obtain optimal solutions for 116 instances, and good upper bounds for 60 instances with a gap of less than 5\% (with an average of 2.8\%). Variations of the algorithm could solve 7 more instances to optimality. 
We provide comprehensive analyses to evaluate the performance of our solution approach.
}

\maketitle




\section{Introduction}

Inefficient distribution networks constitute a huge cost-reduction opportunity for many companies, especially when their distribution network is expanding and poorly managed.
Distribution-related costs considerably affect a company's expenses and constitute about 20\% of total operational costs \citep{SupplyChainChopra}. 
By leveraging their resources, companies become more competitive and more flexible. 
In fact, the optimization of inefficient distribution networks can decrease distribution and transportation costs by up to 25\% \citep{brainandcoreport}.

To optimize distribution networks, companies often employ vendor-managed inventory (VMI) systems.
VMI systems provide information transparency improving coordination and synchronization in the supply chain. The supplier/company controls its customers' inventory level as it has access to the inventory information and demand information of the customers. 
The supplier determines the timing and the amount of delivery to each customer, and the customers' visiting order in its daily vehicle routes. 
As a result, the distribution network can be efficient and well-coordinated, inventory levels are closely controlled, and customer expectations are met. 
More specifically, in VMI systems, the optimization of the distribution network implies making routing and inventory decisions over a planning horizon simultaneously. In the literature \citep{bell1983improving}, such decision-making problems are referred to as inventory-routing problems (IRPs). 

In IRPs, a supplier delivers products to a set of geographically dispersed customers, subject to side constraints, such that the total routing and inventory costs are minimized \citep{coelho2014thirty}. 
In a two-echelon distribution network, intermediate facilities can store and consolidate goods before transferring them to their final customers. This problem is denoted as the two-echelon inventory-routing problem (2E-IRP) and arises in many practical applications, such as city logistics and multi-modal transportation \citep{cuda2015survey}. 
The intermediate facilities, also called satellites \citep{crainic2010two} or distribution centers (DCs), are often located in the city outskirts and act as temporary sources of goods. 
The satellites enable fewer trips between the supplier and the final customers, and allow smaller trucks with lower pollution levels entering the city center. 
The use of satellites is increasingly important as more cities are banning the entry of large trucks to the city center permanently or at certain times.

In this paper, we focus on the 2E-IRP, where multiple suppliers can supply several satellites, which in turn can transfer freight to the final customers to meet their demands in the current and subsequent periods. 
Direct deliveries from the suppliers to the customers are not allowed; however, they can be easily incorporated. 
The 2E-IRP is part of a broader class of two-echelon distribution systems, where the flow of freight in one echelon must be coordinated with that in the other echelon.
This coordination makes the routing problem complex and difficult to be effectively decomposed into two sub-problems (one for each echelon) and  efficiently solved separately \citep{cuda2015survey}.
Hence, advanced algorithms are needed to solve this type of problems, both heuristically or to optimality.

The 2E-VRP is similar to the 2E-IRP, except that the inventory decisions are omitted and the planning horizon (usually) comprises only a single period. 
Several exact algorithms are proposed to optimally solve the 2E-VRP \citep{baldacci2013exact, santos2015branch}, while the state-of-the-art exact method is the branch-price-and-cut algorithm proposed in \cite{marques2020improved}.
Several papers also studied variants of the 2E-VRP, such as the 2E-VRP with Time Windows \citep{dellaert2019branch, Mhamedi2020ABA}, the electric 2E-VRP \citep{breunig2019electric, jie2019two}, and the 2E-VRP with pick-up and delivery \citep{belgin2018two, do2021agile}.
The interested reader is referred to \cite{cuda2015survey} and \cite{sluijk} for extensive reviews on the 2E-VRP.



In this paper, we propose a novel exact branch-and-price algorithm for the 2E-IRP. The contributions are threefold:
\begin{itemize}
\item We propose a mathematical model for the 2E-IRP, where customers can be visited from any satellite. Similar to \cite{desaulniers2016branch}, in which the authors focus on solving the classical IRP exactly, the model does not contain inventory balance constraints for the customers or the satellites. To guarantee flow conservation, we introduce variables detailing the usage of each delivery to the customers and the satellites, proven to yield tighter lower bounds, as in \cite{krarup1977plant}.
\item We develop a branch-and-price algorithm to solve the 2E-IRP. We propose a branching scheme and a labeling algorithm, extended from the algorithm presented in \cite{desaulniers2016branch} to solve new pricing subproblems associated with each satellite and time period.
Each subproblem corresponds to an elementary shortest path problem with resource constraints combined with the linear relaxation of a knapsack problem.
Solving a subproblem returns second-echelon routes and the related deliveries starting and ending at the satellite in the time period associated with the subproblem.
To accelerate the solution method, we solve the pricing subproblems in a multi-threaded fashion and use several acceleration techniques.
\item We generate 400 new instances for the 2E-IRP derived from instances provided in \cite{archetti2007branch}. The original instances involve one supplier, a set of customers varying from 5 to 50, three or six periods, and one vehicle. 
These instances were also used to evaluate multi-vehicle algorithms (up to five vehicles). 
The newly generated instances involve one or two suppliers, two or three satellites, a set of customers varying from five to 25, three periods, and two to five second-echelon vehicles. 
We test our algorithm on these new instances, and provide insights that can be used to improve the performance of our algorithm. 
\end{itemize}


The remainder of the paper is organized as follows. A literature review is given in Section \ref{section LIT REV}. In Section \ref{section Prob DEF}, a formal description and a mathematical model are provided. The proposed branch-and-price algorithm is presented in Section \ref{section BP algorithm}. Section \ref{section Instance} describes the instance design. The results of the computational experiments are discussed in Section \ref{section RESULTS}, which is followed by conclusions in Section \ref{section CONCLUSION}.

\section{Literature Review} \label{section LIT REV}
Many IRP variants are studied over the years (see the comprehensive review by \citealt{coelho2014thirty}), including different inventory policies, among which the most common one is the maximum-level (ML) policy. In this policy, the inventory level is bounded by the capacity available to the customer. 
To solve the IRP variants, a wide range of algorithms are proposed. 
Several papers developed exact methods, where most of them focused on branch-and-cut algorithms \citep{archetti2007branch, solyali2011branch, coelho2013exact, adulyasak2014formulations, Guimares2020MechanismsFF}. Only one paper \citep{desaulniers2016branch} presented a branch-price-and-cut algorithm.
\cite{Guimares2020MechanismsFF} and \cite{desaulniers2016branch} are the state-of-the-art algorithms for the classic IRP.
The former performs best on instances containing one to three vehicles whereas the latter performs best on instances containing four and five vehicles.
\cite{Guimares2020MechanismsFF} embedded two mechanisms in a branch-and-cut algorithm to improve the feasible solutions and render infeasible solutions feasible.
These mechanisms showed a good performance for the IRP and two of its variants.
\cite{desaulniers2016branch} developed an ad-hoc labeling algorithm derived from the one proposed for the split-delivery vehicle routing problem \citep{desaulniers2010branch}. 
The adapted labeling algorithm solves a column generation subproblem per period, which corresponds to an elementary shortest path problem with resource constraints (ESPPRC) combined with the linear relaxation of a knapsack problem.
In this paper, we build on the work of \cite{desaulniers2016branch} and develop a branch-and-price algorithm for the 2E-IRP. 

The literature on the 2E-IRP is rather limited.
The initial papers studying this type of problem, \citep{chan1998probabilistic, zhao2008model, li2011solution}, developed heuristic approaches. 
These papers focused on minimizing the average travel and inventory costs over an infinite planning horizon for a system with one supplier one DC and many customers,
except the work of \cite{chan1998probabilistic} where many DCs are considered.
\cite{chan1998probabilistic} developed an algorithm based on the Fixed Partition (FP) Zero-Inventory Ordering (ZIO) policy.
In the FP policy \citep{anily1993two}, customers are partitioned into disjoint sets and each set of customers is visited by one vehicle. 
In the ZIO policy, a customer receives an order when its inventory level is down to zero.
\cite{zhao2008model} proposed the FP power-of-two (POT) strategy that combines both the FP policy and the POT policy \citep{roundy198598}.
In this strategy, customers are first partitioned based on the FP policy, and then a POT policy is determined, in which each set of customers is restricted to be visited at the same replenishment interval which is power of two times a basic planning period.
They develop a variable large neighborhood search algorithm to optimize the partition sets and solve the 2E-IRP.
\cite{li2011solution} study the same version of 2E-IRP except that direct deliveries from the supplier to the customers are allowed.
The authors developed a decomposition solution approach based on the FP policy. 
They decompose the problem into three subproblems under a given fixed partition and use a genetic algorithm to optimize the partition sets. 
The resulting subproblems are solved separately using efficient heuristics and a dynamic programming algorithm.  

Extensions of the 2E-IRP were recently studied.
For instance, \cite{rohmer2019two} presented a 2E-IRP for perishable products and proposed an adaptive large neighborhood search method to solve it.
\cite{ji2020mixed} introduced a mixed-integer robust programming model for 2E-IRP with perishable products and time windows under demand uncertainty.
\cite{nambirajan2016care} extended a three-phase heuristics CAR (Clustering, Allocation, Routing) presented by \cite{ramkumar2011mathematical} into CAR-Extended for a multi-product 2E-IRP, with multiple suppliers producing mutually exclusive products, one DC, and a set of customers. 
\cite{guimaraes2019two} presented a branch-and-cut algorithm and a matheuristic to solve a variant of 2E-IRP where all decisions are made at the level of the DCs. 
In other words, a homogeneous fleet located at a DC collects the input from one supplier and returns to the DC (with back-and-forth trips only) before delivering the goods to the final customers. 
\cite{schenekemberg2020two} extended the approach proposed by \cite{guimaraes2019two} and designed a branch-and-cut algorithm and a matheuristic to solve the 2E-IRP with fleet management, where short- and mid-term agreements manage the fleet rental. 
\cite{schenekemberg2021two} also studied the two-echelon production-routing problem arising in the petrochemical industry, where production decisions are made at the DCs. 
The authors developed a branch-and-cut algorithm and a novel exact algorithm that employs parallel computing techniques to combine local search procedures within a traditional branch-and-cut scheme.

The above papers focused on variants of the 2E-IRP, and only \cite{farias2021model} addressed a basic version of the 2E-IRP under the VMI setting, and is considered the closest to this paper.
\cite{farias2021model} considered a single supplier (that is also the decision maker), multiple DCs, and a set of customers that are pre-assigned to DCs. 
The authors proposed a branch-and-cut algorithm for the 2E-IRP and compared different routing configurations (one-vehicle, multi-vehicle, and multi-tours). 
They also studied multiple inventory policies including ML policy, order-up-to-level (OU) policy and order fixed quantity policy (OFQ). 
The OU policy assumes that whenever a DC or a customer is replenished, its inventory level must reach its maximum level. 
The OFQ policy imposes that the quantity delivered to every DC (and similarly to every customer) is pre-defined and fixed.
The main differences between the two studies are that they assume a single supplier with a limited inventory capacity and charge inventory holding costs at the supplier.
Moreover, the setting of the instances used in \cite{farias2021model} is also different from the one used in our paper.
The former considers that each subset of customers associated with the same DC constitutes a city, that the cities are distant from each other, and that a single supplier located faraway from the cities supply all customers' demands via the DCs, whereas we consider that the DCs are located at the outskirt of the same city where all customers are located (see Section \ref{section Instance}).
The pre-assignment of customers reduces the complexity of the problem as it limits the second-echelon routes and the quantities entering the DCs. 
In addition to travel costs and inventory holding costs, \cite{farias2021model} also includes fixed ordering costs of the DCs and the customers, that can be easily incorporated in our algorithm.
In the remainder of this paper, we will use the terminology proposed by \cite{crainic2010two} and refer to the DCs by satellites.


\section{Problem Definition and Formulation} \label{section Prob DEF}
Consider a set $U$ of suppliers, a set $S$ of satellites, and a set $N$ of customers. 
Each supplier has an unlimited production capacity over the planning horizon $T = \{1,..., \tau\}$. 
An artificial period $\tau+1$ is considered to handle all end-of-horizon inventories.
Given a homogeneous fleet of $K^1$ first-echelon vehicles of capacity $Q^1$, suppliers transfer freight to a subset of satellites, such that each satellite is visited at most once in period $t\in T$. 
Each satellite $i \in S$ and each customer $c \in N$ have an initial inventory ($I_i^0$ and $I_c^0$), an inventory capacity ($C_i$ and $C_c$), and an inventory holding cost ($f^H_i$ and $f^H_c$).
A homogeneous fleet of $K^2$ vehicles of capacity $Q^2$ is used to deliver freight from satellites to customers, where $Q^2 \leq Q^1$. 
Each customer $c$ is visited at most once in period $t \in T$ and requires a known quantity of freight $d_c^t$.
The objective of this problem is to minimize the travel and inventory holding costs while meeting customers' demands and respecting capacity constraints over the planning horizon. 

Spanning a single period, a first-echelon vehicle travels a first-echelon route $p \in P$ starting from a supplier, visiting a subset of satellites and ending at the same supplier. 
Similarly, a second-echelon route $r \in R$ starts from a satellite, visits a subset of customers and ends at the same satellite. 
No direct delivery from a supplier to the customers is allowed.
We consider the following sequence of operations: at the beginning of each period, first-echelon vehicles can transfer freight from suppliers to satellites. 
The latter see their inventory replenished using an ML inventory replenishment policy. 
Subsequently, second-echelon vehicles pick up the freight from the satellites and deliver it to the customers. 
At the end of the period, the customers consume their demand and holding costs are charged on the remaining inventory. Finally, all vehicles return to their origin and become available for the next period.

The problem can be defined on an undirected graph $G=(V, E)$, where the vertex set $V$ and the edge set $E$ are defined as $V=U \cup S \cup N$ and $E=E^1 \cup E^2$, respectively. 
$E^1$ denotes the set of first-echelon edges and consists of edges linking suppliers to satellites and satellites to each other, $E^1 = \{<i,j> | i\in U \cup S, j\in S \}$. $E^2$ corresponds to the set of second-echelon edges and consists of edges linking satellites to customers and customers to each other, $E^2 = \{<i,j> | i\in S \cup N, j\in N \}$ . The cost of travelling along an edge $e$ is denoted by $f_e$.

Following \cite{desaulniers2016branch}, a path-based formulation is used to model the problem as a mixed-integer program. 
Let $P$ be the set of first-echelon routes, $P_s$ the set of first-echelon routes visiting satellite $s\in S$, and $S_p$ the set of satellites visited by route $p\in P$. 
Similarly, let $R$ be the set of second-echelon routes, $R_s$ the set of second-echelon routes starting at satellite $s \in S$, $R_c$ the set of second-echelon routes visiting customer $c\in N$, and $N_r$ set of customers visited by route $r \in R$. Let $\widetilde{x}_{e}^{p}$ ($\widetilde{x}_{e}^{r}$) be integer parameters indicating the number of times edge $e \in E^1$ ($e\in E^2$) is used in route $p$ ($r$).  

As given in \cite{desaulniers2016branch}, a set $W_r^t$ of route delivery patterns (RDPs) is associated with each route $r$ and each period $t$. 
An RDP $w \in W_r^t$ specifies the deliveries to each customer $c\in N_r$ along route $r$ in period $t$. 
RDP $w$ details the quantities delivered in period $t$ to customer $c$ into sub-deliveries, namely the quantities $(q_{wc}^h)_{t \leq h \leq \tau +1}$ dedicated to cover (fully or partially) the demand of period $h\in T$ or dedicated to remain in the end inventory $\tau+1$.
To limit the number of potential sub-deliveries, it is assumed that, without loss of generality, the consumption of delivered quantities respects the first-in first-out (FIFO) rule. 

Under the FIFO rule, the initial inventory $I_c^0$ of customer $c \in N$ is consumed first to fulfill its demands of the first periods. 
The remaining quantity from the initial inventory at customer $c$ at the end of period $h\in T$ is denoted by $I_c^{0,h} = \max\{0, I_c^0 - \sum_{l=1}^h d_c^l\}$.
The residual demands (the demands not covered by the initial inventory) at customer $c$ are computed as :
\begin{align*}
 \bar{d}_c^h=
    \begin{cases}
    \max\{0, \ d_c^1 - I_c^0\} & \mathrm{if} \ h=1\\
    \max\{0, \ d_c^h - I_c^{0, \; h-1}\}& \mathrm{otherwise, }\ 
    \end{cases}
\quad \quad \quad \quad \quad \forall h\in T.
\end{align*}

Hence, the set of periods associated with the sub-deliveries of a delivery to customer $c$ in period $t$ can be reduced from $\{t, ..., \tau+1\}$ to the set $T_{ct}^+ = \{h\in \{t,t+1, ...., \tau+1\} | (h\in T, \bar{d}_c^h>0 \text{ and } (h=t \text{ or } \sum_{l=t}^{h-1} d_c^l < C_c)) \text{ or } (h=\tau+1 \text{ and } \sum_{l=t}^{h-1} d_c^l < C_c)\}$. 
Moreover, an upper bound $UB_{ct}^h$ on the sub-deliveries $q_{wc}^h$ to customer $c$ for each period $h\in T_{ct}^+$ can be defined as:
\begin{align*}
 \mli{UB}_{tc}^h=
    \begin{cases}
    \min\{\bar{d}_c^h, C_c-I_c^{0,h-1}\},  & \quad \quad \quad \quad \mathrm{if} \ h=t\\
    C_c - \sum_{l=t}^{h-1} d_c^l - I_c^{0,h-1} & \quad \quad \quad \quad \mathrm{if} \ h=\tau +1\\
    \min\{\bar{d}_c^h, C_c - \sum_{l=t}^{h-1} d_c^l -I_c^{0,h-1}\},  & \quad \quad \quad \quad \mathrm{otherwise} .\\
    \end{cases}
\end{align*}
Thus, an RDP $w\in W_r^t$ consists of the sub-deliveries $(q_{wc}^h)_{\substack{c\in N_r\\ h\in T_{ct}^+}}$ where $q_{wc}^h \in [0, UB_{tc}^h]$ and the total quantity delivered in $w$ is $q_w = \sum_{c\in N_r} \sum_{h\in T_{ct}^+} q_{wc}^h$. A sub-delivery $q_{wc}^h$ corresponds to a \textit{zero sub-delivery} if $q_{wc}^h=0$ , a \textit{full subdelivery} if $q_{wc}^h = UB_{ct}^h$ and a \textit{partial sub-delivery} otherwise. As in \cite{desaulniers2016branch}, only\textit{ extreme} RDPs are considered, that is RDPs with at most one partial sub-delivery, and their convex combination is used to generate any other RDP. For this, we associate every route $r$ and RDP $w \in W_r^t$ with a continuous variable $\alpha_{rw}^t$ with values in $[0,1]$ that provides the proportion of the route $r \in R$ operated with RDP $w$.  

Given customer $c \in N$, and period $h\in T$,  three sets of periods essential to our formulation are defined below. Denote by $T_{ch}^- = \{t \in T | h\in T_{ct}^+\}$ the set of periods at which a sub-delivery can be made to fulfill the demand $d_c^h$ of customer $c$ in period $h$. 
Moreover, let $\Gamma_{ch}^- = \{t\in \{1, ...,h\}| \exists k \in \{h, ..., \tau+1\}, k \in T_{ct}^+\}$ denote the set of periods at which a sub-delivery can be made to fulfill the demand of customer $c $ in period $h $ or any subsequent period, including end inventory $\tau+1$. 
Finally, $\Gamma_{ct}^+=\{h\in T | t \in \Gamma_{ch}^-\} $ is the inverse set of  $\Gamma_{ch}^- $. 

The remaining quantity from $q_w$ at the inventory of customer $c$ at the end of period $h \in T_{ct}^+ $ is denoted by $b_{wc}^h = \sum_{\substack{l\geq h+1\\ l\in T_{ct}^+}} q_{wc}^l$. 
It corresponds to the total quantity dedicated to fulfill the demands of customer $c$ in the periods subsequent to $h$.
The cost associated with route $r$ and RDP $w$ can be computed as $f_{rw} = \sum_{e\in E^2} f_e\widetilde{x}_{e}^{r} + \sum_{c\in N_r} \sum_{t \leq h\leq T } f_c^H b_{wc}^h$, where the first term corresponds to the travel costs and the second term to the inventory holding costs at customers in route $r$ and related to RDP $w$. 
The binary requirements on second-echelon routes are not imposed directly on $\alpha_{rw}^t$, but rather on additional variables $\alpha_r^t$, that represent the usage of a second-echelon route $r$ in period $t$.

Binary variable $\lambda_p^t$ is equal to 1 if first-echelon route $p\in P$ is selected in period $t$ and 0 otherwise. 
The cost associated with first-echelon route $p\in P$ is computed as $f_p = \sum_{e\in E^1} f_e \widetilde{x}_{e}^{p}$. 
As in the facility location-based formulation introduced by \cite{krarup1977plant}, transfer quantities at the satellites are modeled as non-negative variables $\psi_{s,t}^l$ that indicate the quantity entering satellite $s$ in period $l\in \{0\} \cup T$ and exiting it in period $t \in T$ or remaining in the end inventory in period $t = \tau +1$ such that $l\leq t$. 
The initial inventory at the satellite $s$ can be formulated as : $I_s^0 = \sum_{t=1}^{\tau+1} \psi_{s,t}^0 $. 

The mathematical formulation of the 2E-IRP is written as follows:

\begin{alignat}{2}
&\min  \sum_{t \in T}\left[\sum_{p \in P}  f_{p} \ \lambda_{p}^{t} + \sum_{r \in R} \sum_{w \in W_r^t} \ f_{rw} \ \alpha_{r, w}^{t} + \sum_{s\in S} \sum_{l=0}^{t} \sum_{h=t+1}^{\tau+1} f_s^H \psi_{s,h}^l \right] \span  \label{MP.1}\\ \intertext{s.t.}
&\sum_{l=0}^{t}\psi_{s,t}^l = \sum_{r \in R_s} \sum_{w \in W_r^t}   q_w \: \alpha_{r,w}^t & \forall s \in S, t \in T \label{MP.2}\\
& \sum_{t\in T_{ch}^-} \sum_{r\in R_c} \sum_{w\in W_r^t} q_{wc}^h \alpha_{rw}^t = \bar{d}_{c}^h  &  \forall c \in N,  h \in T \text{ such that  } \bar{d}_c^h >0 \label{MP.3} \\
& \sum_{l=0}^{t} \sum_{k= t}^{\tau+1}\psi_{s,k}^l \leq C_{s}   &\quad \quad   \forall s \in S, t \in T \label{MP.4} \\
&I_{c}^{0, h} + \sum_{t \in \Gamma_{ch}^-}\sum_{r \in R_{c}} \sum_{w \in W_r^t} \sum_{\substack{l \in T_{ct}^+\\ l>h}}q_{w, c}^l \; \alpha_{r,w}^t + d_{c}^h\leq C_{c}   &\quad \quad   \forall c \in N,  h \in T \label{MP.5} \\
&\sum_{s \in S_{p}} \sum_{k= l}^{\tau+1}\psi_{s,k}^l \leq Q^{1} + Q^1 (|S_p| -1 )(1-\lambda_p^l) & \quad \quad   \forall p \in P, \  l \in T \label{MP.6} \\
&\sum_{p \in P_s} \lambda_{p}^{t}  \leq 1 & \quad \quad  \forall s \in S, t \in T  \label{MP.7} \\
&\sum_{r \in R_{c}}  \sum_{w \in W_r^t} \alpha_{r, w}^{t}  \leq 1 & \quad \quad   \forall c \in N,\ \forall t \in T \label{MP.8} \\
&\sum_{p \in P} \lambda_{p}^{t}  \leq K^{1} & \quad    \forall t \in T  \label{MP.9} \\
&\sum_{r \in R}  \sum_{w \in W_r^t} \alpha_{r, w}^{t}\leq K^{2} & \quad    \forall t \in T  \label{MP.10} \\
&\sum_{t= 1}^{\tau+1} \psi_{s,t}^0 = I_s^0
& \quad \quad \forall s\in S  \label{MP.11} \\
&\sum_{t= l}^{\tau+1} \psi_{s,t}^l \leq Q^1 \sum_{p\in P_s} \lambda_p^l & \quad \quad \forall s\in S, l \in T  \label{MP.12} \\
&\sum_{p\in P_s} \lambda_p^l \leq  \sum_{t= l}^{\tau+1} \psi_{s,t}^l  & \quad \quad \forall s\in S, l \in T  \label{MP.19} \\
&\alpha_r^t= \sum_{w\in W_r^t} \alpha_{r,w}^{t}  & \quad    \forall r \in R, \  t \in T \label{MP.15}\\
&\psi_{s,t}^l \geq 0 & \forall s\in S,\ t \in T\cup\{\tau+1\},\ 0 \leq l \leq t\label{MP.13} \\
&\alpha_{r,w}^{t} \geq 0 & \quad    \forall r \in R, \  t \in T,  \ w \in W_r^t \label{MP.14} \\
&\alpha_r^t\in \{0, 1\}  & \quad    \forall r \in R, \  t \in T  \label{MP.16}\\
&\lambda_{p}^{t} \in \{0,1\} & \quad   \forall p \in P, \  t \in T,  \label{MP.17}
\end{alignat}

The objective function (\ref{MP.1}) minimizes the sum of transportation costs of first-echelon routes (first term), total costs of second-echelon routes including the inventory holding costs at the customers (second term), and the inventory holding costs at the satellite (third term). 
Constraints (\ref{MP.2}) define the flow of goods exiting satellite $s$ in period $t$, and link variables $\alpha_{rw}^t$ to $\psi_{s,t}^l$. Constraints (\ref{MP.3}) ensure that the demand of each customer is met in each period.
Constraints (\ref{MP.4}) and (\ref{MP.5}) impose that the inventory capacities of the satellites and the customers are not exceeded (ML policy).
Constraints (\ref{MP.6}) ensure that first-echelon vehicles capacity is not exceeded.
Constraints (\ref{MP.7}) and (\ref{MP.8}) enforce that satellites and customers are visited at most once during a period. 
Constraints (\ref{MP.9}) and (\ref{MP.10}) ensure that the number vehicles used in each echelon does not exceed the number of available vehicles. 
Constraints (\ref{MP.11}) set the initial inventory of the satellites.
Constraints (\ref{MP.12}) ensure that no quantity enters to a satellite if it is not visited by a first-echelon route.
Similarly, constraints (\ref{MP.19}) ensure that a satellite receives at least one item if it is visited by a first-echelon route.
Constraints (\ref{MP.15}) define the binary variables $\alpha_r^t$ that indicates whether route $r$ is selected.
Constraints (\ref{MP.13})-(\ref{MP.17}) define the domain of the decision variables. 

The above formulation (\ref{MP.1})-(\ref{MP.17}) involves a large number of variables $\alpha_{rw}^t$ associated with second-echelon routes. 
The second-echelon routes and the associated RDPs are dynamically generated in a column generation procedure described in Section \ref{Column Generation}.
On the other hand, the number of first-echelon variables $\lambda_p^t$ depends on the number of suppliers, satellites, and periods considered. 
The problem setting employed in this paper makes these numbers take moderate values, hence the enumeration of all first-echelon routes is possible.
Moreover, as suppliers have an unlimited production capacity, the set of first-echelon routes can be reduced to the non-dominated first-echelon routes by solving the traveling salesman problem for each subset $\mathcal{S} \subset S$ and each supplier $u\in U$ (as in \citealt{marques2020improved}), and selecting the shortest route starting from any supplier and visiting all satellites in the subset $\mathcal{S}$. 
The resulting number of variables $\lambda_p^t$ is equal to $|T|(2^{|S|}-1)$. 
This procedure can be executed as a pre-processing step before the start of the solution method elaborated in the next section. 

\section{Branch-and-Price Algorithm} \label{section BP algorithm}
We propose a branch-and-price algorithm to solve the 2E-IRP. 
At each node of the search tree, a linear programming relaxation of problem (\ref{MP.1})-(\ref{MP.17}), also called the master problem (MP), is solved using a column generation approach.
This section describes the column generation approach and the branching strategy. 

\subsection{Column Generation}\label{Column Generation}
To solve the linear relaxation of problem (\ref{MP.1})-(\ref{MP.17}), constraints (\ref{MP.15}) are not needed and are omitted. 
The MP consists, then, of the model (\ref{MP.1})-(\ref{MP.19}), (\ref{MP.13}) and (\ref{MP.14}). 

In the first iteration of the column generation procedure, the set $P$ of first-echelon routes is initialized with all non-dominated routes, and the set $R$ of second-echelon routes is initialized with a small set of feasible routes $r$ and their associated set $W_r^t$ of RDPs. 
These initial solutions are generated using a greedy heuristic that, for period $t \in T$, selects a satellite $s$ and generates round-trips between $s$ and each customer $c \in N$ to deliver its residual demand $\bar{d}_c^t$ at period $t$, until the total quantity exiting the satellite $s$ is about to exceed its capacity $C_s$ or the capacity $Q^1$ of a first-echelon vehicle. 
Another satellite is, then, selected and the same procedure is repeated until all customers are visited in each period $t \in T$. 
Artificial variables with a very large cost are added to the model to ensure feasibility since the initial routes might not yield a feasible solution to the MP. 

The MP restricted to the set of initial solutions, the restricted master problem (RMP), is then solved and its dual values are used to define the pricing subproblems.
Upon solving the subproblems, either new columns with negative reduced costs are found and added to the RMP, triggering a new iteration of the column generation procedure, or no such column is found, and thus optimality to the linear problem is proven and the column generation is stopped. 

The pricing subproblems can be decomposed into $|S|*|T|$ subproblems, where a subproblem $SP_s^t$ is defined for each satellite $s\in S$ and each time period $t \in T$. 
In the next subsections, we define subproblem $SP_{s}^t$, and the labelling algorithm applied to solve it.

\subsubsection*{Subproblem Definition.}
Given a satellite $s \in S$ and a time period $t \in T$, let us define subproblem $SP_s^t$. Subproblem $SP_s^t$ minimizes the reduced cost of columns, where a column corresponds to a feasible route $r\in R_s$ starting and ending at satellite $s$ and a feasible extreme RDP $w \in W_r^t$. 
A route $r \in R_s$ is considered feasible if it visits each customer at most once, and the total quantity delivered along route $r$ does not exceed the second-echelon vehicle's capacity $Q^2$. 

Consider the dual variables denoted by $\pi_{st}^{\ref{MP.2}}$, $\pi_{ch}^{\ref{MP.3}}$, $\pi_{ch}^{\ref{MP.5}}$, $\pi_{ct}^{\ref{MP.8}},$ and $ \pi_{t}^{\ref{MP.10}}$ associated with constraints (\ref{MP.2})-(\ref{MP.3}), (\ref{MP.5}), (\ref{MP.8}), and (\ref{MP.10}), respectively. The reduced cost $\bar{f}_{rw}^t$ of $\alpha_{rw}^t$ is computed as follows:
\begin{align}
\bar{f}_{rw}^t = f_{rw}^t +q_w \pi_{st}^{\ref{MP.2}} - \sum_{c\in N_r} \sum_{\substack{h \in T_{ct}^+\\ h\neq \tau+1}} q_{wc}^h \pi_{ch}^{\ref{MP.3}}- \sum_{c\in N_r}\sum_{h \in \Gamma_{ct}^+} \sum_{\substack{l\in T_{ct}^+\\ l>h}} q_{wc}^l\ \pi_{ch}^{\ref{MP.5}} - \sum_{c\in N_r} \pi_{ct}^{\ref{MP.8}} - \pi_{t}^{\ref{MP.10}}
\end{align}

Subproblem $SP_s^t$ can be represented on the graph denoted by $(V^t_s, A_s^t)$, where the set of vertices $V^t_s=N\cup \{s^{src}, s^{snk}\}$ contains the set of customers $N$, a source vertex $s^{src}$ and a sink vertex $s^{snk}$ representing the satellite $s$ at the start and at the end of period $t$, respectively. 
The set of arcs $A_s^t=\{(i,j)|i\in N\cup\{s^{src}\}, j\in N\cup \{s^{snk}\}\}$ contains the arcs linking the customers to each other and the satellite $s$ to all customers. The reduced cost $\bar{f}_{ij}$ of arc $(i,j) \in A_s^t$ is : 
\begin{align}
\bar{f}_{ij}=
    \begin{cases}
    f_{ij} - \pi_t^{(\ref{MP.10})}  & \quad \quad \mathrm{if }\ i= s^{src} \\
    f_{ij} - \pi_{it}^{(\ref{MP.8})} & \quad \quad \mathrm{otherwise,} \label{5.29}
\end{cases}
\end{align}
To formulate the objective function of the subproblem, two variables are required, one for the routing part and the other for the delivery part. 
Let $x_{ij}^r$ be a binary variable indicating whether the arc $(i,j)\in A_s^t$ is selected in route $r$ and let $\xi_c^h$ be an integer variable with a value in $[0, UB_{ct}^h]$, representing the quantity delivered to customer $c \in N$ in period $t$ and dedicated to period $h \in T_{ct}^+$. 
The total quantity delivered along route $r \in R_s$ should not exceed the second-echelon vehicle's capacity $Q^2$, and the quantity $\xi_c^h$ should be equal to 0 if the customer $c$ is not visited in route $r$.
The objective function of subproblem $SP_s^t$ can be written as:
\begin{align}
\bar{f}_s^{t}(x, \xi) &= \sum_{(i,j) \in A_s^t} \bar{f}_{i,j} \: x_{i,j} + \sum_{c\in N_r} \sum_{h\in T_{ct}^+} \xi_c^h  (\pi_{st}^{\ref{MP.2}} -\pi_{ch}^{\ref{MP.3}} - \sum_{\substack{l \in \Gamma_{ct}^+ \\ l<h} } \pi_{cl}^{\ref{MP.5}} + \sum_{\substack{t\leq l\leq T\\ l< h}} f_c^H) \label{eq:21}
\end{align}
where $\pi_{ch}^{\ref{MP.3}}$ is 0 if $h=\tau +1$. 

The subproblem $SP_s^t$ can be seen as an ESPPRC combined with the linear relaxation of a knapsack problem (see \citealp{desaulniers2010branch} for a similar definition).
Labeling algorithms used to solve classic vehicle routing problems and their variants cannot be applied to solve this subproblem because the delivery quantity to each customer is a decision variable ($\xi_c^h$), yielding a reduced cost and a load resource which are functions of these decision variables. 

\subsubsection*{Labeling Algorithm and Accelerating Techniques.}
The labeling algorithm proposed in \cite{desaulniers2016branch} is extended to solve the subproblem $SP_s^t$, where the reduced cost of the customer delivery patterns (CDPs) are derived from Equation (\ref{eq:21}) (see Appendix \ref{Appendix A} for the full definition).
A CDP specifies a combination of subdelivery types associated with customer $i\in N$ in period $t \in T$.
It can be seen as a component of the sequence representing an extreme RDP.
Note also that some of the dominance rules were adapted as well with respect to the reduced cost of a label (see Appendix \ref{Appendix A}).
The following accelerating techniques are used to improve the performance of the solution method along with customer delivery patterns handling (or CDP handling) as described in \cite{desaulniers2016branch}.
The CDP handling is performed at the beginning of each column generation iteration and for each subproblem, and the dominance rules used to compare CDPs follow the adapted dominance rules.
\begin{itemize}
  \item \textit{Bidirectional dynamic programming:}
\end{itemize}  
\cite{righini2006symmetry} proposed bidirectional labeling that consists of propagating labels in both directions (forward and backward) until reaching half of the critical resource, and joining forward and backward labels associated with the same vertex. 
In our paper, the considered critical resource is the second-echelon vehicle capacity, $Q^2$.
\cite{righini2008new} provide a detailed explanation of the implementation of the bidirectional labeling. 
\begin{itemize}
  \item \textit{Symmetry break:}
\end{itemize}
The labelling algorithm generates symmetric routes (routes that can be traversed in both directions). 
To avoid generating symmetric routes and later discarding one of them by using dominance rules, we impose that the index of the first customer visited along a route is less than that of its last customer as in \cite{pessoa2009robust}.
\begin{itemize}
  \item \textit{Ng-path relaxation:}
\end{itemize}
The \textit{ng-path} relaxation accelerates the pricing problem solving and is used in state-of-art branch-price-and-cut algorithms \citep{costa2019exact}.
\cite{baldacci2011new} proposed the \textit{ng-path} relaxation which defines a neighborhood $N_c$ for each customer $c\in N_s^t$ that contains $\kappa$ closest customers to $c$, including $c$ itself, in terms of distance ($\kappa$ is set to 5 in our computational tests).
An \textit{ng-path} can contain a cycle starting and ending at customer $c^{'}$ if and only if there exists a customer $c$ in this cycle such that $c^{'} \notin N_c$.

\begin{itemize}
  \item \textit{Parallel computing:}
\end{itemize}
Solving the pricing problems is the most time consuming task in column generation. 
For the instances used in this paper, the number of subproblems solved ranges from nine to 30 in each iteration of the column generation procedure. 
It is beneficial to solve them in parallel using multi-threading as the subproblems are independent of each other and the time spent to solve them can compensate for the time required to create multiple threads.
\begin{itemize}
  \item \textit{Heuristic column generators:}
\end{itemize}
Solving the pricing problem using the exact labelling algorithm can be very time-consuming, especially in the first iterations of the column generation and when artificial variables are active. 
Heuristic column generators can accelerate the column generation procedure drastically.
Three heuristic column generators are used, where one method is based on the tabu search as described in \cite{archetti2011enhanced}, and the other two are heuristic labeling algorithms.
The first heuristic labeling algorithm consists of solving the labelling algorithm for a subset of customers.
The subset of customers is chosen such that the reduced cost of the roundtrip routes $(s^{src}$- $c$ - $s^{snk})$ is among the $k$ least reduced costs ($k$ is equal to 5 in our case). 
The second heuristic labeling algorithm corresponds to the labeling algorithm applied on a reduced network. 
A label associated with customer $c_1$ can be extended to only one customer $c_2$ such that the resulting label has the least reduced cost.
To compute the extension, the CDP with the lowest cost is selected assuming no quantity is delivered in the partial deliveries.

The pseudocode for the column generation procedure at a given branching node is provided in Algorithm \ref{alg:1}.  
Initialization is performed (lines \ref{alg line 2}-\ref{alg line 6}).
In a given iteration of the column generation procedure, the subproblems are all solved either heuristically or exactly, and the resulting columns are added to the RMP. Subsequently, the RMP is optimized and a new iteration begins.
This is implemented in two steps.
The first step (lines \ref{alg line 7}-\ref{alg line 9}) consists of solving all subproblems heuristically as many times as needed until no subproblem returns a solution.
In the second step (lines \ref{alg line 10}-\ref{alg line 28}), subproblems are solved alternately, once exactly then heuristically, until the stopping criterion is met, i.e. the model is no longer improvable. 
When solving exactly (line \ref{alg line 14}), if no appropriate solution is found for all subproblems, then the column generation stops (line \ref{alg line 16}). 
Otherwise, the solutions with negative reduced costs are added to the RMP and the RMP is optimized again (line \ref{alg line 28}).
If no solution is found for one or more subproblems, these subproblems are no longer solved heuristically or exactly (sp.solve\_exact = False).
The remaining subproblems are then alternately solved heuristically and exactly (lines \ref{alg line 18}-\ref{alg line 27}) until all subproblems have not yielded a solution in an exact solve iteration (lines \ref{alg line 24}-\ref{alg line 25}).
In that case, all subproblems are again allowed to be solved (lines \ref{alg line 12}-\ref{alg line 13}), and the second step is repeated until all subproblems yield no solution in a single iteration.
When solving heuristically, a subproblem is solved only if it is allowed to be solved exactly (line \ref{alg line 32}).
In that case, it is first solved using tabu search (line \ref{alg line 33}). 
If tabu search does not yield any solution with a negative reduced cost, the first heuristic labeling algorithm is applied in the same iteration (lines \ref{alg line 34}-\ref{alg line 35}). 
If no appropriate solution is found (line \ref{alg line 36}), the first heuristic labeling algorithm
is no longer applied in the subsequent iterations (line \ref{alg line 38}), and the second heuristic labeling algorithm 

{\fontsize{11}{15}\selectfont
\begin{algorithm}[H]

\DontPrintSemicolon
\SetAlgoLined
\caption{Column Generation}\label{alg:1}
\KwData{SPs: list of subproblems} 
\SetKwFunction{FMain}{Column\_Generation}
\SetKwProg{Fn}{Procedure}{:}{\KwRet}
\Fn{\FMain{}}{
solve\_h := False; last\_iter := True \label{alg line 2} \; 
solve\_e $:=$ True; model\_improvable $:=$ True\;
solutions $:= \emptyset $ \label{alg line 4} \;
\For{\upshape sp $ \in $SPs }{  \label{alg line 5}
    sp.solve\_exact $:=$ True, sp.solve\_h1 $:=$ True} \label{alg line 6}
 \While (\tcp*[f]{First step}) {\upshape solutions $\neq \emptyset $}{ \label{alg line 7}
  solutions := H\_solve() \;
  Add solutions with negative reduced costs to the RMP and solve the new LP 
  } \label{alg line 9}
 \While (\tcp*[f]{Second step}){ \upshape model\_improvable = True}{\label{alg line 10}
    \uIf{ \upshape solve\_e = True and last\_iter = True}{ \label{alg line 11}
        \For{\upshape sp $\in$ SPs}{\label{alg line 12}
            sp.solve\_exact := True} \label{alg line 13}
        solutions := E\_solve() \label{alg line 14} \tcp*{Multithreading solve}
        \eIf{\upshape solutions $= \emptyset$}{\label{alg line 15}
            model\_improvable := False \label{alg line 16}}
            {last\_iter := False; solve\_e := False; solve\_h := True}\label{alg line 18}
    }
    \uElseIf{ \upshape solve\_h = True \label{alg line 19}}{
        solutions := H\_solve() \label{alg line 20} \;
        solve\_e := True; solve\_h := False} \label{alg line 21}
    \Else{
        solutions := E\_solve() \tcp*{Multithreading solve}\label{alg line 23}
        \eIf{ \upshape solutions $= \emptyset$}{ \label{alg line 24}
            last\_iter := True \label{alg line 25}}{
            solve\_e := False; solve\_h := True}}\label{alg line 27}
Add solutions with negative reduced costs to the RMP and solve the new LP \label{alg line 28}
 }}
\BlankLine
\SetKwFunction{Hsolve}{H\_solve}
\SetKwProg{Pn}{Function}{:}{\KwRet}
\Pn{\Hsolve{}}{
solutions $:= \emptyset$\;
\For{ \upshape sp $\in$ SPs}{
    \If{ \upshape sp.solve\_exact = True}{ \label{alg line 32}
    sol := solutions returned by tabu search \label{alg line 33}\;
    \uIf{ \upshape sol $= \emptyset$ and  sp.solve\_h1 = True}{\label{alg line 34}
        sol := solutions returned by the heuristic 1 \label{alg line 35}\;
        \If{ \upshape sol $= \emptyset$}{ \label{alg line 36}
            sol := solutions returned by the heuristic 2 \label{alg line 37}\;
            sp.solve\_h1 := False \label{alg line 38}}
            }
    \uElseIf{ \upshape sol $= \emptyset$ and  sp.solve\_h1 = False}{ \label{alg line 39}
    sol := solutions returned by the heuristic 2 \label{alg line 40}\;}
    solutions := solutions  $\cup$ sol \label{alg line 41}\;
}}
}\Return solutions
\end{algorithm}
}

\noindent  is applied in the same iteration (line \ref{alg line 37}).
The routes with negative reduced costs are added to the RMP and a new iteration of the column generation starts.




\subsection{Branching}
Branching can be required on first-echelon routes, on second-echelon routes, or on both. Consequently, 10 types of branching decisions are imposed, namely:
\begin{enumerate}
    \item The total number of first-echelon routes over all periods ($\sum_{t\in T} \sum_{p\in P} \lambda_p^t$). \label{BD.1}
    \item The number of first-echelon routes in each period $t\in T$ ($\sum_{p\in P} \lambda_p^t$).  \label{BD.2}
    \item The flow through each satellite $s \in S$ from first-echelon routes over all periods ($\sum_{t \in T}\sum_{p\in P_s} \lambda_p^t$). \label{BD.3}
    \item The use of first-echelon routes $p \in P$ in each period $t \in T$ ($\lambda_p^t$).\label{BD.4}
    \item The total number of second-echelon routes over all periods ($\sum_{t\in T} \sum_{r\in R}\sum_{w \in W_r^t} \alpha_{r,w}^t$). \label{BD.5}
    \item The total number of second-echelon routes in each period $t\in T$ ($\sum_{r\in R}\sum_{w \in W_r^t} \alpha_{r,w}^t$). \label{BD.6}
    \item The flow through each customer $c\in C$ over all periods ($\sum_{t\in T} \sum_{r\in R_c}\sum_{w\in W_r^t} \alpha_{r,w}^t$).  \label{BD.7}
    \item The flow through each customer $c\in C$ in each period $t \in T$ ($\sum_{r\in R_c}\sum_{w\in W_r^t} \alpha_{r,w}^t$).  \label{BD.8}
    \item The flow through each customer $c\in C$ in each period $t \in T$ from each satellite $s\in S$ ($\sum_{r\in R_c\cap R_s}\sum_{w\in W_r^t} \alpha_{r,w}^t$).\label{BD.9}
    \item The flow on each edge $<i, j>$ of the second-echelon in each period $t\in T$ ($\sum_{r \in R} \sum_{w \in W_r^t} (x_{ij}^r + x_{ji}^r) \ \alpha_{r,w}^t$). \label{BD.10}
\end{enumerate}
The first four branching decisions ensure the integrality requirements (\ref{MP.17}) on variables $\lambda_p^t$, and can be imposed by adding a constraint in the RMP. 
Such addition does not require any change in the subproblems. 
The remaining branching decisions enforce the integrality requirements (\ref{MP.16}) on variables $\alpha_{rw}^t$.
Similarly, they can be imposed by adding a constraint to the RMP, however the dual values of these constraints must be incorporated in the reduced costs of certain arcs. 
Moreover, in the last branching decision, both arcs $(i,j)$ and $(j,i)$ must be removed from $A_s^t$ when the flow on an edge $<i,j>$ must be set to 0.

When the solution of the linear relaxation solution is fractional, we first check whether branching decisions on first-echelon routes (\ref{BD.1} to \ref{BD.4}) can be applied.
In that case, we first branch on decision \ref{BD.1} when possible, otherwise we compute the value of each candidate variable and branch on the one whose fractional value is closest to 0.5. 
If none of the branching decisions on first-echelon routes can be applied, we compute the value of the candidate variables for the remaining decisions (\ref{BD.5} to \ref{BD.10}) and choose the candidate for each type whose fractional value is closest to 0.5.
If one of the variables selected for decisions \ref{BD.7}, \ref{BD.8}, \ref{BD.9} or \ref{BD.10} has a fractional value in the interval $[0.25, 0.75]$ then we branch on one of these variables with priority given to \ref{BD.7}, \ref{BD.8}, \ref{BD.9} and then \ref{BD.10}. Otherwise, we branch on the candidate variable whose fractional value is closest to 0.5. 

The search tree is explored using a best-first strategy. 
In this strategy, the node selection consists of always choosing the most promising node, i.e., the node with the least lower bound, among the list of open nodes of the search tree.

\section{Instance design} \label{section Instance}
The IRP instances of \cite{archetti2007branch} involve one supplier and a number of customers ranging from 5 to 50 for a planning horizon of 3 periods, and a number of customers ranging from 5 to 30 for a planning horizon of 6 periods.
\cite{archetti2007branch} created 160 instances, divided into four classes based on the inventory holding costs (high H and low L) and on the length of the planning horizon (3 and 6): H3, H6, L3, L6.

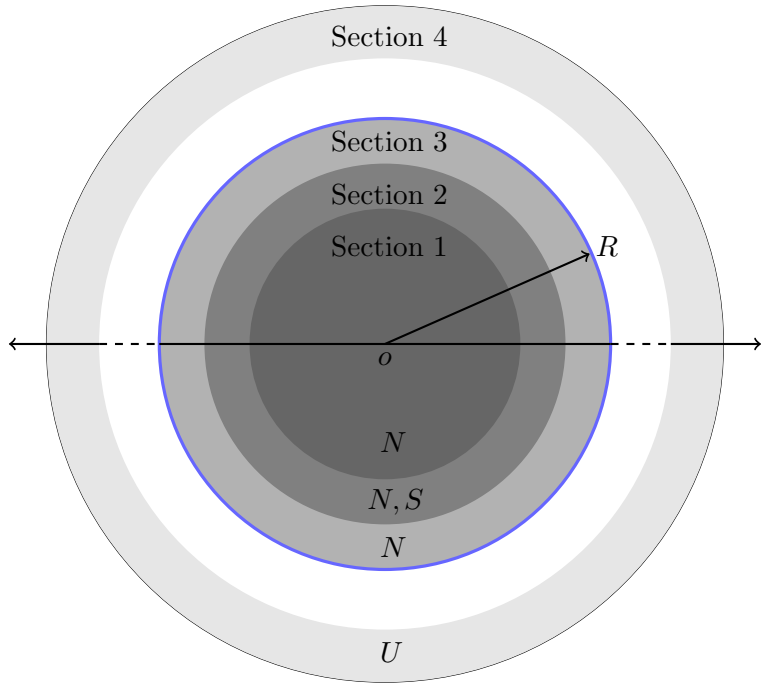
\begin{figure}[!h] 
  \caption{The schematic representation of an urban area, adapted from \cite{dellaert2019branch} }
\label{fig urban area}
\centering 
\begin{tikzpicture}
\draw[thick,<->] (0,4) -- (8,4);
\draw (4,4) circle (1.8cm);
\draw (4,4) circle (2.4cm);
\draw (4,4) circle (3cm);
\draw (4,4) circle (3.8cm);
\draw (4,4) circle (4.5cm);
\draw[color=white] (4,4) circle (5.1cm);
\fill[black!10!white] (4,4) circle (4.5cm);
\fill[color=black,fill=white, very thick] (4,4) circle (3.8cm);
\filldraw[color=blue!60,fill=black!30!white, very thick] (4,4) circle (3cm);
\fill[black!50!white] (4,4) circle (2.4cm);
\fill[black!60!white] (4,4) circle (1.8cm);
\draw[thick,<-] (-1,4) -- (0.2,4);
\draw[thick,dashed, -] (0.2,4) -- (1,4);
\draw[thick,-] (1,4) -- (7,4);
\draw[thick,dashed, -] (7,4) -- (7.8,4);
\draw[thick,->] (7.8,4) -- (9,4);
\draw[thick,->] (4,4) -- (6.72,5.2);
\node[right] at (3.15,5.3) {Section 1};
\node[right] at (3.15,6) {Section 2};
\node[right] at (3.15,6.7) {Section 3};
\node[right] at (3.15,8.1) {Section 4};
\node[right] at (3.8,2.7) {$N$};
\node[right] at (3.63,1.9) {$N, S$};
\node[right] at (3.8,1.3) {$N$};
\node[right] at (3.8,-0.1) {$U$};
\node[] at (4,3.8)  {$o$};
\node[] at (6.95,5.3)  {$R$};
\end{tikzpicture}
\end{figure}

To derive 2E-IRP instances from the original IRP instances of \cite{archetti2007branch}, a transformation is performed. 
To distinguish between the original instances' and our instances' properties, we use the same notation as described in Section \ref{section Prob DEF} and add a superscript $A$ to refer to Archetti's instances.
We simulate a circular urban area that is divided into four sections (see Figure \ref{fig urban area}). 
First, the outer limit of the urban section 3 (the blue circle) is defined as the smallest circle containing all customers. 
Its center $o$ represents the city center and its radius $R$ represents the width of the housing area, i.e., all customers are located within this circle.
Urban section 2 is the area hosting the satellites, and represents the area surrounding the city center.
It is simulated by a ring whose center is $o$ and whose smaller and larger radii are 90\% and 99\% of $R$, respectively.
This section is divided into $|S|$ subsections and each satellite is randomly located in one of these subsections.

\noindent Similarly, section 4 is the area hosting suppliers and represents the external area of the city.
It is simulated by a ring whose center is $o$ and whose smaller and larger radii are 250\% and 300\% of $R$, respectively.
The external area is divided into $|U|$ subsections such that each supplier is located in one of these subsections.
A similar procedure was implemented to generate instances for the 2E-VRP with Time Windows in \cite{dellaert2019branch}. 

We consider the classes with a planning horizon of three periods and the instances with up to 25 customers.
For each class, we create two combinations \textit{a}s\textit{b} (1s2 and 2s3) based on the number of suppliers $a \in \{1,2\}$ and the number of satellites $b \in \{2,3\}$. 
We set the satellites' capacity as the sum of the supplier's capacity $C_{u}^A$ and the quantity made available at the supplier $u^A$ in each period. 
The satellites' initial inventory is proportional to the total residual demand of all customers over all periods divided by the number of satellites.
For each satellite residual demand, a coefficient of proportionality is randomly chosen between 0.4 and 0.6.
Moreover, the number of the first-echelon trucks $K^1$ is equal to the number of the suppliers $|U|$, and their capacity is twice of the quantity made available at the supplier $u^A$. 
The number of second-echelon trucks $K^2$ ranges from 2 to 5, and their capacity is simply the original vehicle capacity divided by the number of trucks $K^2$.
The resulting new data set is composed of 400 instances in total.

\section{Numerical Results} \label{section RESULTS}
The branch-and-price algorithm is tested on instances derived from \cite{archetti2007branch} as explained above. In this section, we provide an analysis of the results for the newly generated instances. All tests are performed on 32 cores of AMD Rome 7H12 processors clocked at 2.6 GHz and with 2 GB of memory per core.
The solution method is implemented in Python 3.7.
CPLEX v20.1 is used to solve the RMPs. 
We first discuss the results in terms of integrality gaps, integer solutions, and then provide insights for each class of instances. 


\subsection{Integrality Gaps}

We solve the newly generated instances by implementing our branch-and-price algorithm with a three-hour time limit. 
We solve the integer RMP with CPLEX after solving the root node and after solving 20 nodes when at least 20 nodes are solved within the time limit and are needed for the algorithm to terminate.
The reason for solving the integer RMP at the root node and after 20 nodes is that the algorithm often does not find an upper bound at the termination in our preliminary experiments. 
The results are presented by number $K^2$ of second-echelon vehicles and by groups of instances, where a group contains instances of the same class with the same combination of suppliers and satellites (1s2 or 2s3).
\begin{table}[!h]
\caption{Average Gaps (in Percentage) per Instance Group and Number of Vehicles}
\label{tab:tab 1}
\resizebox{\textwidth}{!}{%
\begin{tabular}{@{}llrrrrrrrrrrrrrrrrr@{}}
\toprule
\multicolumn{2}{l}{\multirow{2}{*}{}} &
  \multicolumn{5}{c}{\multirow{2}{*}{Gap$_0$}} &
  \multirow{2}{*}{} &
  \multicolumn{5}{c}{\multirow{2}{*}{Gap$_{20}$}} &
  \multirow{2}{*}{} &
  \multicolumn{5}{c}{\multirow{2}{*}{Gap$_f$}} \\
\multicolumn{2}{l}{} & \multicolumn{5}{l}{}                  &  & \multicolumn{5}{l}{}                  &  & \multicolumn{5}{l}{}                \\ \cmidrule(lr){3-7} \cmidrule(lr){9-13} \cmidrule(l){15-19} 
Instance \\  Class &
  Combination &
  $K^2$=2 &
  $K^2$=3 &
  $K^2$=4 &
  $K^2$=5 &
  Average &
   &
  $K^2$=2 &
  $K^2$=3 &
  $K^2$=4 &
  $K^2$=5 &
  Average &
   &
  $K^2$=2 &
  $K^2$=3 &
  $K^2$=4 &
  $K^2$=5 &
  Average \\ \midrule
H3      & 1s2        & 57.70  & 51.12 & 47.11 & 44.04 & 49.84 &  & 7.71  & 7.65  & 6.37  & 5.64  & 6.72  &  & 10.05 & 6.31  & 4.71 & 3.89 & 6.16  \\
        & 2s3        & 57.78 & 54.27 & 50.29 & 47.60  & 52.38 &  & 12.37 & 11.96 & 11.96 & 11.91 & 12.02 &  & 14.22 & 10.19 & 7.68 & 4.99 & 9.17  \\
        & Average H3 & 57.74 & 52.69 & 48.70  & 45.82 & 51.11 &  & 9.96  & 9.86  & 9.30   & 8.90   & 9.43  &  & 12.14 & 8.25  & 6.19 & 4.44 & 7.67  \\
L3      & 1s2        & 63.27 & 55.18 & 50.51 & 45.58 & 53.33 &  & 6.85  & 7.24  & 6.91  & 6.05  & 6.71  &  & 10.42 & 6.21  & 4.85 & 3.70  & 6.17  \\
        & 2s3        & 59.97 & 60.39 & 56.78 & 53.52 & 57.60  &  & 11.95 & 13.96 & 15.16 & 15.24 & 14.33 &  & 12.25 & 11.84 & 9.94 & 7.18 & 10.24 \\
        & Average L3 & 61.62 & 57.79 & 53.65 & 49.55 & 55.47 &  & 9.31  & 10.60  & 11.23 & 10.64 & 10.55 &  & 11.33 & 9.03  & 7.40  & 5.44 & 8.20   \\
Average &            & 59.64 & 55.24 & 51.17 & 47.69 & 53.27 &  & 9.64  & 10.25 & 10.26 & 9.79  & 10.00    &  & 11.74 & 8.64  & 6.80  & 4.94 & 7.93  \\ \bottomrule
\end{tabular}%
}
\end{table}

Table \ref{tab:tab 1} presents the average integrality gaps for each instance group with respect to the number of vehicles. 
This table includes the average integrality gap at the root node ($\text{Gap}_0$), the average integrality gap after solving 20 nodes ($\text{Gap}_{20}$), and the final integrality gap ($\text{Gap}_f$).
The integrality gap $\text{Gap}_0$ ($\text{Gap}_{20}$; $\text{Gap}_f$) is computed as ($\bar{z} - \underbar{$z$})/ \underbar{$z$}$, where $\bar{z}$ is the lower bound at the root node (after solving 20 nodes; at the termination of the algorithm) and $\underbar{$z$}$ is the upper bound at the root node (after solving 20 nodes; at the termination of the algorithm).

The integrality gap at the root node is relatively high, on average 53.27\% with a minimum of 25\% and a maximum of 109\%. 
The solution of the linear relaxation at the root node includes fractional values of the variables $(\lambda_p^t)_{\substack{p\in P\\ t\in T }}$ associated with first-echelon routes, which have high transportation costs. 
The average absolute improvement of the integrality gap after solving 20 nodes, which is computed as ($\text{Gap}_0-\text{Gap}_{20}$) when possible, is 42.56\%, with a minimum of 23.76\% and a maximum of 81.62\%. 
In our preliminary tests, we find that a sharp decrease occurs at the children nodes of the root node mainly by means of the branching decisions on first-echelon routes. 
The branching decision of type \ref{BD.1} is the first branching decision to be applied after the root node is solved and can be seen as a valid inequality on the minimum number of first-echelon routes. 
The integrality gap decreases to 7.93\% on average at the end of the three-hour limit.

We observe that the lower the number of vehicles is the higher the average integrality gap, except for $\text{Gap}_{20}$. 
However, this exception can be explained by the fact that for most difficult instances (instances with a high number of customers and 2 to 3 vehicles), we could not solve up to 20 nodes in the search tree within the time limit. 
The reported results for $\text{Gap}_{20}$ correspond to relatively easy instances, hence are biased. 
We also observe that the combination with one supplier and two satellites (1s2) has lower average integrality gaps ($\text{Gap}_0, \text{Gap}_{20}$, and $\text{Gap}_f$) than the combination with two suppliers and three satellites (2s3). 
This result was expected because when the number of suppliers and satellites increases, the number of possible solutions (and variables) increases, and the lower bound is less tight. Finally, we observe slightly higher average integrality gaps for L3 than H3.

\subsection{Integer Solution Results}
We analyze the results by integrality gap ranging from 0 to 0.05\%, 0.05\% to 5\%, greater than 5\%, and when no feasible solution was found. 
We assume that an instance was solved to optimality if the integrality gap is less than 0.05\%.  
Table \ref{tab:tab 2} shows that we could solve up to 116 instances to optimality.
In addition, 60 instances were solved with a gap smaller than 5\% (an average of 2.8\%) and are considered as good feasible solutions. 
We were able to find a feasible solution for 214 instances that have a gap higher than 5\% and an average of 13.67\%.
Finally, we could not find a feasible solution for 10 instances, all of which have 25 customers and two vehicles.
For these instances, the time limit was insufficient even to solve the root node.
This observation is consistent with the results obtained for the one-echelon inventory-routing problem \citep{desaulniers2016branch}, 
\begin{table}[h!]
\caption{Number of Instances by Integrality Gap Range and Instance Group}
\label{tab:tab 2}
\begin{tabular}{@{}llrrrr@{}}
\toprule
Instance Class & Combination & Optimal Solution & Gap$_f <5$\% & Gap$_f\geq5$\% & No Solution \\ \midrule
H3    & 1s2         & 30               & 15        & 53         & 2                      \\
      & 2s3         & 27               & 13        & 58         & 2                      \\
      & H3 Subtotal    & 57               & 28        & 111        & 4                      \\
L3    & 1s2         & 31               & 23        & 43         & 3                      \\
      & 2s3         & 28               & 9         & 60         & 3                      \\
      & L3 Subtotal    & 59               & 32        & 103        & 6                      \\
Total &             & 116              & 60        & 214        & 10                     \\ \bottomrule
\end{tabular}
\end{table}
\noindent where the poor performance of their exact solution approach on the two-vehicle  instances is mainly due to the large time spent on solving the subproblems.
We also observe that 1s2 instances are easier to solve in general, and the performance of our algorithm does not depend on the magnitude of the holding costs.

Next, we compare the results in terms of the number of instances solved to optimality by instance group and number of vehicles (see Table \ref{tab:tab 3}).
We solve 116 instances out of the 400 instances generated, most of which have 5 (80 instances) to 10 customers (33 instances). 
From these results, we observe that a general trend with respect to number of vehicles cannot be obtained.
We solve more instances with combination 1s2 than with 2s3 and more instances in L3 than in H3.
Moreover, two factors that indicate the difficulty of an instance are the size of the search tree and the difficulty of the subproblems, which increases as the number of vehicles decreases.
The latter implies an increase in the second-echelon vehicle's capacity leading to a higher number of labels, resulting in larger times to solve the subproblems.

\begin{table}[h!]
\centering
\caption{Number of Instances Solved to Optimality by Instance Group and Number of Vehicles}
\label{tab:tab 3}
\begin{tabular}{@{}llrrrrr@{}}
\toprule
Instance class & Combination & $K^2$=2 & $K^2$=3 & $K^2$=4 & $K^2$=5 & Total \\ \midrule
H3             & 1s2         & 7     & 9     & 6     & 8     & 30    \\
               & 2s3         & 6     & 8     & 6     & 7     & 27    \\
               & H3 Subtotal    & 13    & 17    & 12    & 15    & 57    \\
L3             & 1s2         & 8     & 8     & 7     & 8     & 31    \\
               & 2s3         & 7     & 8     & 6     & 7     & 28    \\
               & L3 Subtotal    & 15    & 16    & 13    & 15    & 59    \\
Total &    & 28    & 33    & 25    & 30    & 116   \\ \bottomrule
\end{tabular}
\end{table}

\begin{table}[]
\caption{Average Computational Times (in Seconds) and Number of Nodes by Instance Group and Gap Range}
\label{tab:tab 4}
\resizebox{\textwidth}{!}{%
\begin{tabular}{@{}llrrrlrrrlrrr@{}}
\toprule
               &             & \multicolumn{3}{c}{Optimal Solution} &  & \multicolumn{3}{c}{Gap$_f<5$\%} &  & \multicolumn{3}{c}{Gap$_f \geq 5$\%} \\ \cmidrule(lr){3-5} \cmidrule(lr){7-9} \cmidrule(l){11-13} 
Instance Class & Combination & \multicolumn{1}{l}{\#Nodes}    & \multicolumn{1}{l}{Time$_{root}$}    & \multicolumn{1}{l}{Time}    &  & \multicolumn{1}{l}{\#Nodes}  & \multicolumn{1}{l}{Time$_{root}$}  & \multicolumn{1}{l}{Time} &  & \multicolumn{1}{l}{\#Nodes}    & \multicolumn{1}{l}{Time$_{root}$}   & \multicolumn{1}{l}{Time}   \\ \midrule
H3    & 1s2        & 381.73  & 8.13  & 1,667.58 &  & 614.73  & 224.88 & 10,800 &  & 196.7  & 1,498.16 & 10,800 \\
      & 2s3        & 569.44  & 6.82  & 945.18  &  & 1,054.38 & 167.12 & 10,800 &  & 177.76 & 1,399.99 & 10,800 \\
      & H3 Average & 470.65  & 7.51  & 1,325.39 &  & 818.86  & 198.06 & 10,800 &  & 186.8  & 1,446.86 & 10,800 \\
L3    & 1s2        & 570.48  & 39.48 & 1,408.18 &  & 722.48  & 160.27 & 10,800 &  & 133.42 & 1,576.42 & 10,800 \\
      & 2s3        & 1,778.21 & 6.56  & 1,682.73 &  & 1,184.11 & 75.65  & 10,800 &  & 298.77 & 1,305.38 & 10,800 \\
      & L3 Average & 1,143.64 & 23.86 & 1,538.47 &  & 852.31  & 136.47 & 10,800 &  & 229.74 & 1,418.54 & 10,800 \\
Total &            & 812.95  & 15.82 & 1,433.77 &  & 836.7   & 165.21 & 10,800 &  & 207.47 & 1,433.23 & 10,800 \\ \bottomrule
\end{tabular}}
\end{table}
In Table \ref{tab:tab 4}, we report the average computational times ($\text{Time}_{root}$ and Time) and the average number of nodes (\#Nodes) by instance group and gap ranges. 
These results show that the number of nodes in the search tree increases in the instances with combination 2s3, especially for L3 instances where it triples compared to the number of nodes in the instances with combination 1s2. 
The average times to solve the root node ($\text{Time}_{root}$) show an increase of magnitude by at least a factor of 4 (most cases by a factor of 10) between instances solved to optimality and solved with $\text{Gap}_{f} < 5\% $ and by at least a factor of 6 for instances solved with $\text{Gap}_{f} < 5\% $ and $\text{Gap}_{f} \geq 5\% $. 
From these results, it is clear that reducing the solving times of the root nodes and more concretely of the subproblems is necessary to solve larger instances. 
The average times to find optimal solutions are, in general, higher for L3 instances, even though 1s2 instances in the H3 group require a significant amount of time. 
This exception is observed mainly due to the larger instances in L3 class that can be solved.
The results for same-size instances show that the average number of nodes, total times, and times at the root node tend to be higher in general for 2s3 instances. 
This observation can be explained by the increase in the number of subproblems (6 for 1s2 instances and 9 for 2s3 instances) and the increase of first-echelon routes. 
We also observe that the $\text{Time}_{root}$ decreases on average when the number of vehicles decreases, an indication that the branch-and-price could be further improved to focus on solving instances with a large number of vehicles.
\begin{table}[h!]
\centering
\caption{Number of Instances and Average Gaps (in Percentage) by Integrality Gap Range and Instance Group for Best Results}
\label{tab:tab 5}
\resizebox{\textwidth}{!}{%
\begin{tabular}{@{}llrlrrlrrlr@{}}
\toprule
 &
   &
  \multicolumn{1}{c}{Optimal Solution} &
  \multicolumn{1}{c}{} &
  \multicolumn{2}{c}{Gap$_f<5$\%} &
  \multicolumn{1}{c}{} &
  \multicolumn{2}{c}{Gap$_f\geq5$\%} &
  \multicolumn{1}{c}{} &
  \multicolumn{1}{c}{No Solution} \\ \cmidrule(lr){3-3} \cmidrule(lr){5-6} \cmidrule(lr){8-9} \cmidrule(l){11-11} 
Instance Class &
  Combination &
  \multicolumn{1}{l}{\#Instances} &
   &
  \multicolumn{1}{l}{\#Instances} &
  \multicolumn{1}{l}{Gap$_f$} &
   &
  \multicolumn{1}{l}{\#Instances} &
  \multicolumn{1}{l}{Gap$_f$} &
   &
  \multicolumn{1}{l}{\#Instances} \\ \midrule
H3    & 1s2      & 31  &  & 16 & 2.80 &  & 51  & 10.58 &  & 4  \\
      & 2s3      & 28  &  & 14 & 2.76 &  & 56  & 14.91 &  & 2  \\
      & H3 Total & 59  &  & 30 & 2.78 &  & 107 & 12.84 &  & 2  \\
L3    & 1s2      & 34  &  & 21 & 3.32 &  & 42  & 11.85 &  & 6  \\
      & 2s3      & 30  &  & 9  & 2.66 &  & 58  & 16.04 &  & 3  \\
      & L3 Total & 64  &  & 30 & 3.12 &  & 100 & 14.28 &  & 3  \\
Total &          & 123 &  & 60 & 2.95 &  & 207 & 13.54 &  & 10 \\ \bottomrule
\end{tabular}%
}
\end{table}

Finally, multiple experiments are conducted with different sets of parameter values to investigate their effect on the results. 
We use the following sets of parameters values: $\kappa \in [5,7,8]$, half-point $\in \{0.5,0.05\}$, best-first search or \textit{local depth-first} search as described in \cite{desaulniers2016branch}.
For each instance, we provide the best results in terms of $\text{Gap}_f$ among all these experiments in Appendix \ref{Appendix B}. 
An overview of these results is presented in Table \ref{tab:tab 5}.
We do not observe a significant difference in the number of instances solved to optimality for each experiment. 
However, we are able to find a few more optimal solutions as we introduce 123 optimal solutions to the literature.
Finally, we observe that \textit{ng-path} relaxation helps solve subproblems for larger instances and is essential to our algorithm. 
Without \textit{ng-path} relaxation, the root node is not solved for 173 instances within a 3-hour time limit.

\section{Conclusions}\label{section CONCLUSION}

We studied a two-echelon multi-depot inventory-routing problem that considers multiple suppliers and customers not pre-allocated to satellites. 
This type of problem often arises in city logistics. 
For this problem, we proposed a route-based formulation and developed a branch-and-price algorithm for solving it. 
For each satellite and each period, a pricing subproblem is defined and is solved by a labeling algorithm resulting in second-echelon routes starting and ending at the same satellite in the same period, whereas first-echelon routes are generated by enumeration. 
Finally, we derived several branching rules for this problem.

We conducted experiments on newly generated instances derived from benchmark instances taken from the literature. 
We introduced 123 optimal solutions and an upper bound for 267 instances with an average gap of 11.16\%.
For the one-echelon inventory-routing problem, the existing exact methods can consistently solve instances with up to 25 customers and three periods within reasonable times.
For the two-echelon case, the problem complexity increases, and more computational effort is required to solve larger instances.

Our work can be extended to include direct deliveries, where customers can be directly visited from the suppliers when they are located in the outskirts of the city. 
Other possible extensions can be introducing location decisions, time windows, or a collect option especially when the customers require, for sustainability reasons, maintenance services for reusable products, which can be done at the supplier or the satellites (e.g., closed-loop supply chains). 
Another extension can focus on accelerating the solution method to solve larger instances by deriving valid inequalities specific to the two-echelon case, or improving the proposed formulation. 
Finally, heuristics and matheuristics approaches based on FP policies can be developed to solve large-size instances.

\ACKNOWLEDGMENT{%
This publication is part of the project Sustainable Supply Chain Management in Healthcare (SSCMH) (with project number 439.18.457) which is financed by the Dutch Research Council (NWO).
This work made use of the Dutch national e-infrastructure with the support of the SURF Cooperative using grant no. EINF-2033.
}

%
%
%
\newpage
\begin{APPENDICES}

\section{Labeling Algorithm} \label{Appendix A}

Given a subproblem $SP_s^t$, let $r$ be a partial path in $(V^t_s, A_s^t)$ from $s^{src}$ to a vertex $i\in N_s^t$ and $w\in W_r^t$ an associated RDP.
We use the terminology proposed by \cite{desaulniers2016branch} for the component of the labeling algorithm. 
Hence, the label denoted $E_i = (T_i^{cost}, T_i^{loadF}, (T_i^{custk})_{k\in N}, T_i^{part}, \\ T_i^{ratePiP}, T_i^{maxP})$ representing the feasible path $r$ with its associated RDP $w$ contains the following components:

\setlength{\tabcolsep}{0.3em}
\begin{tabular}{ r l }
T$_i^{cost}$: & The reduced cost of the path/RDP ($r, w$).\\
T$_i^{loadF}$: & The total quantity delivered along the path $r$ according to RDP $w$.\\
T$_i^{custk}$: & A binary value that indicates whether a customer $k\in N$ has been visited along path $r$.\\
T$_i^{part}$: & A binary value that indicates whether RDP $w$ contains a partial subdelivery.\\
T$_i^{ratePiP}$: & The unit rate of contribution to reduced cost associated with the partial subdelivery if any.\\
T$_i^{maxP}$: & The maximum quantity that can be delivered in the partial subdelivery if any.
\end{tabular}
\\

Similarly, we define for each CDP $\gamma \in \Gamma_{it}$ and each period $h\in T_{it}^+$ a binary parameter $f_{\gamma}^h$ (respectively, $g_{\gamma}^h$) that takes value 1 if CDP $\gamma$ contains a full (respectively, partial) subdelivery for period $h$ and we associate the following values:

\setlength{\tabcolsep}{0.3em}
\begin{tabular}{r l}
$\tau_{\gamma}^{cost}=$ & $\sum_{h\in T_{it}^+} f_{\gamma}^h \ \mli{UB}_{ti}^h \  (\pi_{st}^{\ref{MP.2}} -\pi_{ih}^{\ref{MP.3}} - \sum_{\substack{l \in \Gamma_{it}^+ \\ l<h} } \pi_{il}^{\ref{MP.5}} + \sum_{\substack{t\leq l\leq T\\ l< h}} f_i^H)$  : \text{The contribution to the path/RDP}\\ & reduced cost (\ref{eq:21}) \\
$\tau_{\gamma}^{loadF}= $& $\sum_{h\in T_{it}^+} f^h_{\gamma} \mli{UB}_{ti}^h$ : \text{The total quantity delivered in the full subdeliveries.}\\
$\tau_{\gamma}^{part}=$ & $\sum_{h\in T_{it}^+} g^h_{\gamma}$ : \text{The number of partial deliveries (0 or 1)}.\\
$\tau_{\gamma}^{ratePiP}=$ & $\sum_{h\in T_{it}^+} g_{\gamma}^h   (\pi_{st}^{\ref{MP.2}} -\pi_{ih}^{\ref{MP.3}} - \sum_{\substack{l \in \Gamma_{it}^+ \\ l<h} } \pi_{il}^{\ref{MP.5}} + \sum_{\substack{t\leq l\leq T\\ l< h}} f_i^H)$  : \text{The rate of contribution to the  path/RDP} \\& reduced cost (\ref{eq:21}) for each unit delivered in the partial delivery if any. \\
$\tau_{\gamma}^{maxP}=$ &$ \sum_{h\in T_{it}^+} g^h_{\gamma} (\mli{UB}_{ti}^h -1)$ : \text{The maximum quantity that can be delivered in the partial delivery}\\ &if any.
\end{tabular}
\\

Any CDP $\gamma$ that contains a partial subdelivery such that $\tau_{\gamma}^{ratePiP} \geq 0$ can be discarded as replacing it with a zero subdelivery will always yield in a CDP that has a contribution to the reduced cost at least as good as the zero subdelivery. We assume that the  $\tau_{\gamma}^{ratePiP} < 0$ in the following.

We keep the same extension functions as in \cite{desaulniers2016branch}, however a few changes were made to the dominance rules.
We rewrite the reduced cost $\bar{f}_i$ associated with label $E_i$ as a function of the quantity delivered in the partial subdelivery $\xi^{part}$:
\[\bar{f}_i(\xi^{part}) =  T_i^{cost} \ + \xi^{part}\  T_i^{ratePiP}, \quad \quad \quad \forall \xi^{part} \in [0, T_i^{maxP}]\]
The reduced cost of label $E_i$ corresponds to a line segment and the dominance rules must allow the comparison of line segments.
We use the dominance rule that was introduced by \cite{desaulniers2010branch}, which allows the comparison of two line segments (cost functions associated with labels $E_1$ and $E_2$).
To declare that a label $E_1$ dominates label $E_2$, one of the sufficient conditions related to the cost comparison is that the cost of every feasible path obtained by extending $E_2$ is greater than or equal to the cost of the path obtained by extending $E_1$ similarly, which is satisfied if :
\begin{align*}
    &\bar{f}_1(T_1^{maxP}) = \min_{\xi^{part}\in [0, T_1^{maxP}]} \bar{f}_1(T_1^{\xi^{part}}) \leq \min_{\xi^{part}\in [0, T_2^{maxP}]} \bar{f}_2(T_2^{\xi^{part}}) = \bar{f}_2(T_2^{maxP})\\
    &\bar{f}_1(T_1^{\xi^{part}}) \leq \bar{f}_2(T_2^{\xi^{part}}), \quad \quad \forall \xi^{part} \in [0, T_2^{maxP}]
\end{align*}

The first condition can be interpreted as a comparison of the best costs that can be yielded by both labels, see condition \ref{cond d} below, whereas the second condition is a comparison of the label costs over all  values that can be taken by label $E_2$, see conditions \ref{cond e}-\ref{cond f}.

Thus, the dominance rules can be defined as follows:
\begin{definition}
A label $E_1 = (T_1^{cost}, \ T_1^{loadF}, \ (T_1^{custk})_{k\in N}, \ T_1^{part}, \ T_1^{ratePiP}, \ T_1^{maxP})$ is said to dominate a label $E_2 = (T_2^{cost}, \ T_2^{loadF}, \ (T_2^{custk})_{k\in N}, T_2^{part}, \ T_2^{ratePiP}, \ T_2^{maxP})$ if both labels are associated with the same vertex and the following conditions are satisfied:
\begin{enumerate}[label=(\alph*)]
\item $T_1^{loadF} \leq T_2^{loadF}; $
\item $T_1^{custk} \leq T_2^{custk}, \quad \forall k\in N;$
\item $T_1^{part} \leq T_2^{part}; $
\item $T_1^{cost} + T_1^{maxP} \ T_1^{ratePiP} \leq T_2^{cost} + T_2^{maxP} \ T_2^{ratePiP}; $ \label{cond d}
\item $T_1^{cost} \leq T_2^{cost}; $ \label{cond e}
\item $T_1^{cost} + T_2^{maxP} \ T_1^{ratePiP} \leq T_2^{cost} + T_2^{maxP} \ T_2^{ratePiP}. $ \label{cond f}
\end{enumerate}
\end{definition}

Dominated labels according to this dominance rule can be discarded except when both labels dominate each other, in which case one of the two labels must be kept.
\newpage
\section{Results}  \label{Appendix B}
\footnotesize  
\setlength\LTleft{-30pt}            
\setlength\LTright{-30pt} 
\begin{longtable}{@{\extracolsep{\fill}}lrrrrrrrr@{}}
\caption{Results for H3 Instances with One Supplier and Two Satellites}
\label{tab:H3 1s2}\\
\toprule
Instance &
  \multicolumn{1}{l}{Time} &
  \multicolumn{1}{l}{Gap$_f$} &
  \multicolumn{1}{l}{\#Nodes} &
  \multicolumn{1}{l}{LB} &
  \multicolumn{1}{l}{UB} &
  \multicolumn{1}{l}{Time$_{root}$} &
  \multicolumn{1}{l}{Gap$_0$} &
  \multicolumn{1}{l}{Gap$_{20}$} \\* \midrule
\endfirsthead
\endhead
\bottomrule
\endfoot
\endlastfoot
1s2\_1n5\_2k  & 27.05     & $0^{}$      & 105   & 2,806.39 & 2,806.39 & 1.25      & 56.94 & 4     \\
1s2\_1n5\_3k  & 57.01     & $0^{}$      & 299   & 3,123.32 & 3,123.32 & 1         & 48.51 & 5.37  \\
1s2\_1n5\_4k  & 44.65     & $0^{}$      & 269   & 3,306.97 & 3,306.97 & 0.81      & 35.48 & 2.50  \\
1s2\_1n5\_5k  & 95.23     & $0^{}$      & 689   & 3,549.29 & 3,549.29 & 0.47      & 33.62 & 3.57  \\
1s2\_1n10\_2k & 10,800    & $1.48^{1}$  & 509   & 3,362.48 & 3,412.20 & 46.34     & 53.66 & 5.71  \\
1s2\_1n10\_3k & 9,439.39  & $0^{}$      & 1,393 & 3,618.75 & 3,618.75 & 18.15     & 46.11 & 7     \\
1s2\_1n10\_4k & 10,800    & $1.43^{1}$  & 1,495 & 3,839.77 & 3,894.86 & 12.61     & 44.71 & 3.73  \\
1s2\_1n10\_5k & 10,800    & $1.46^{}$   & 2,969 & 4,054.14 & 4,113.46 & 9.27      & 43.04 & 4.11  \\
1s2\_1n15\_2k & 10,800    & $7.14^{}$   & 55    & 4,091.33 & 4,383.42 & 525.27    & 53.18 & 8.74  \\
1s2\_1n15\_3k & 10,800    & $5.1^{1}$   & 135   & 4,559.51 & 4,792.23 & 233.46    & 42.67 & 6.48  \\
1s2\_1n15\_4k & 10,800    & $5.56^{1}$  & 145   & 4,966.52 & 5,242.61 & 247.08    & 39.57 & 7.12  \\
1s2\_1n15\_5k & 10,800    & $8.22^{1}$  & 99    & 5,384.22 & 5,827.01 & 196.12    & 41.26 & 9.87  \\
1s2\_1n20\_2k & 10,800    & $13.24^{}$  & 7     & 4,020.59 & 4,553.01 & 2,520.97  & 56.93 & N/A   \\
1s2\_1n20\_3k & 10,800    & $5.36^{}$   & 15    & 4,361.47 & 4,595.41 & 1,493.82  & 45.29 & N/A   \\
1s2\_1n20\_4k & 10,800    & $5.63^{}$   & 81    & 4,668.52 & 4,931.37 & 478.15    & 43.88 & 6.45  \\
1s2\_1n20\_5k & 10,800    & $1.85^{3}$  & 41    & 4,980.44 & 5,072.81 & 400.88    & 37.73 & 3.14  \\
1s2\_1n25\_2k & 10,800    & $17.71^{}$  & 3     & 4,423.11 & 5,206.28 & 5,566.13  & 54.59 & N/A   \\
1s2\_1n25\_3k & 10,800    & $17.87^{3}$ & 5     & 4,805.63 & 5,664.28 & 3,913.24  & 53.24 & N/A   \\
1s2\_1n25\_4k & 10,800    & $13.22^{3}$ & 7     & 5,132.84 & 5,811.37 & 2,704.69  & 47.92 & N/A   \\
1s2\_1n25\_5k & 10,800    & $5.22^{}$   & 29    & 5,610.90 & 5,903.61 & 922.88    & 38.72 & 5.35  \\
1s2\_2n5\_2k  & 65.14     & $0^{}$      & 335   & 2,290.40 & 2,290.40 & 1.01      & 75.82 & 8     \\
1s2\_2n5\_3k  & 26.24     & $0^{}$      & 161   & 2,359.54 & 2,359.54 & 0.71      & 67.37 & 3.74  \\
1s2\_2n5\_4k  & 24.40     & $0^{}$      & 173   & 2,527.50 & 2,527.50 & 0.56      & 72.42 & 4.16  \\
1s2\_2n5\_5k  & 11.17     & $0^{}$      & 83    & 2,756.96 & 2,756.96 & 0.45      & 73.78 & 1.22  \\
1s2\_2n10\_2k & 10,800    & $5.93^{3}$  & 525   & 4,130.89 & 4,375.70 & 51.97     & 70.17 & 12.43 \\
1s2\_2n10\_3k & 6,342.61  & $0^{}$      & 427   & 4,568.21 & 4,568.21 & 15.56     & 56.46 & 7.76  \\
1s2\_2n10\_4k & 10,800    & $5.22^{}$   & 1,915 & 4,953.33 & 5,212    & 13.95     & 63.14 & 11.53 \\
1s2\_2n10\_5k & 5,022.66  & $0^{}$      & 815   & 5,265.34 & 5,265.34 & 11.16     & 47.79 & 5.01  \\
1s2\_2n15\_2k & 10,800    & $11.5^{}$   & 47    & 4,176.70 & 4,656.91 & 557.10    & 62.02 & 12.68 \\
1s2\_2n15\_3k & 10,800    & $7.14^{1}$  & 35    & 4,560.65 & 4,886.25 & 564.73    & 49.57 & 7.78  \\
1s2\_2n15\_4k & 10,800    & $5.96^{}$   & 93    & 5,174.47 & 5,482.91 & 98.94     & 46.76 & 7.73  \\
1s2\_2n15\_5k & 10,800    & $4.12^{}$   & 183   & 5,568.82 & 5,798.31 & 89.52     & 40.14 & 6.05  \\
1s2\_2n20\_2k & 10,800    & $15.38^{}$  & 5     & 4,428.26 & 5,109.26 & 3,990.10  & 62.09 & N/A   \\
1s2\_2n20\_3k & 10,800    & $13.72^{3}$ & 11    & 4,886.64 & 5,557.26 & 1,421.97  & 59.53 & N/A   \\
1s2\_2n20\_4k & 10,800    & $6.06^{}$   & 21    & 5,241.18 & 5,558.86 & 774.25    & 47.15 & 6.20  \\
1s2\_2n20\_5k & 10,800    & $5.41^{}$   & 35    & 5,678.73 & 5,985.70 & 468.59    & 47.01 & 5.91  \\
1s2\_2n25\_2k & 10,800    & N/A     & 1     & N/A      & N/A      & 10,810.61 & N/A   & N/A   \\
1s2\_2n25\_3k & 10,800    & $16.74^{}$  & 5     & 4,904.70 & 5,725.95 & 2,972.06  & 52.77 & N/A   \\
1s2\_2n25\_4k & 10,800    & $10.73^{}$  & 15    & 5,590.96 & 6,190.98 & 1,093.72  & 53.98 & N/A   \\
1s2\_2n25\_5k & 10,800    & $5.33^{}$   & 17    & 6,042.19 & 6,363.98 & 706.43    & 45.68 & N/A   \\
1s2\_3n5\_2k  & 5.41      & $0^{}$      & 23    & 3,300.81 & 3,300.81 & 0.75      & 60.17 & 3.19  \\
1s2\_3n5\_3k  & 42.67     & $0^{}$      & 337   & 3,831.81 & 3,831.81 & 0.52      & 77.41 & 10.42 \\
1s2\_3n5\_4k  & 182.31    & $0^{}$      & 1,531 & 4,262.80 & 4,262.80 & 0.39      & 73.94 & 7.69  \\
1s2\_3n5\_5k  & 6.79      & $0^{}$      & 59    & 4,299.20 & 4,299.20 & 0.43      & 65.04 & 1.20  \\
1s2\_3n10\_2k & 10,800    & $6.95^{3}$  & 229   & 3,543.41 & 3,789.82 & 185.87    & 63.58 & 9.75  \\
1s2\_3n10\_3k & 10,800    & $3.25^{1}$  & 331   & 3,983.88 & 4,113.33 & 35.47     & 54.15 & 6.51  \\
1s2\_3n10\_4k & 10,800    & $5.9^{}$    & 901   & 4,410.21 & 4,670.45 & 16.17     & 53.06 & 9.28  \\
1s2\_3n10\_5k & 10,800    & $5.58^{1}$  & 1,259 & 4,812.54 & 5,081.18 & 8.10      & 50.77 & 7.41  \\
1s2\_3n15\_2k & 10,800    & $6.58^{3}$  & 29    & 4,485.36 & 4,780.56 & 1,266.16  & 51.88 & 7.61  \\
1s2\_3n15\_3k & 10,800    & $2.67^{1}$  & 59    & 4,924.06 & 5,055.29 & 417.81    & 45.01 & 4.09  \\
1s2\_3n15\_4k & 10,800    & $2.55^{}$   & 213   & 5,211.31 & 5,344.16 & 98.94     & 44.50 & 5.14  \\
1s2\_3n15\_5k & 6,987.97  & $0^{}$      & 137   & 5,543.70 & 5,543.70 & 54.39     & 37.31 & 1.28  \\
1s2\_3n20\_2k & 10,800    & $19.81^{}$  & 11    & 4,252.22 & 5,094.62 & 1,646.75  & 62.29 & N/A   \\
1s2\_3n20\_3k & 10,800    & $11.26^{}$  & 31    & 4,696.41 & 5,225.38 & 800.79    & 48.26 & 11.39 \\
1s2\_3n20\_4k & 10,800    & $6.18^{1}$  & 19    & 5,174.11 & 5,493.98 & 899.32    & 37.65 & 6.22  \\
1s2\_3n20\_5k & 10,800    & $5.67^{3}$  & 33    & 5,589.50 & 5,906.19 & 664.88    & 35.07 & 6.12  \\
1s2\_3n25\_2k & 10,800    & $56.96^{}$  & 1     & 3,406.02 & 5,345.98 & 10,218.79 & 56.96 & N/A   \\
1s2\_3n25\_3k & 10,800    & $14.49^{}$  & 5     & 4,778.27 & 5,470.52 & 2,750.22  & 51.41 & N/A   \\
1s2\_3n25\_4k & 10,800    & $11.08^{}$  & 9     & 5,202.05 & 5,778.34 & 1,599.23  & 51.84 & N/A   \\
1s2\_3n25\_5k & 10,800    & $6.22^{}$   & 19    & 5,599.60 & 5,947.66 & 1,002.59  & 48.28 & 6.25  \\
1s2\_4n5\_2k  & 58.49     & $0^{}$      & 145   & 2,702.75 & 2,702.75 & 1.29      & 56.49 & 9.95  \\
1s2\_4n5\_3k  & 27.66     & $0^{}$      & 133   & 3,035.09 & 3,035.09 & 0.69      & 47.74 & 8.24  \\
1s2\_4n5\_4k  & 71.49     & $0^{}$      & 427   & 3,402.61 & 3,402.61 & 0.78      & 37.44 & 5.71  \\
1s2\_4n5\_5k  & 21.67      & $0^{}$      & 139    & 3,820.17 & 3,820.17 & 0.33      & 31.4 & 2.81  \\
1s2\_4n10\_2k & 10,268.75 & $0^{}$      & 469   & 4,115.18 & 4,115.18 & 49.73     & 49.06 & 6.10  \\
1s2\_4n10\_3k & 4,242.49  & $0^{}$      & 501   & 4,547.57 & 4,547.57 & 21.83     & 45.53 & 12.28 \\
1s2\_4n10\_4k & 3,895.15  & $0^{}$      & 607   & 4,959.61 & 4,959.61 & 10.30     & 44.20 & 4.92  \\
1s2\_4n10\_5k & 10,800    & $5.69^{}$   & 2,889 & 5,427.28 & 5,736.08 & 7.46      & 49.27 & 12.09 \\
1s2\_4n15\_2k & 10,800    & $12.92^{}$  & 29    & 3,729.01 & 4,210.87 & 563.44    & 63.56 & 13.35 \\
1s2\_4n15\_3k & 10,800    & $3.68^{}$   & 49    & 4,037.65 & 4,186.43 & 511.46    & 44.13 & 5.21  \\
1s2\_4n15\_4k & 10,800    & $4.14^{}$   & 75    & 4,318.97 & 4,497.67 & 194.28    & 39.93 & 5.44  \\
1s2\_4n15\_5k & 10,800    & $2.9^{1}$   & 53    & 4,593.04 & 4,726.33 & 196.92    & 32.72 & 3.54  \\
1s2\_4n20\_2k & 10,800    & $14.55^{1}$ & 5     & 4,759.74 & 5,452.08 & 4,692.03  & 53.10 & N/A   \\
1s2\_4n20\_3k & 10,800    & $14.17^{}$  & 19    & 5,241.64 & 5,984.21 & 774.13    & 52.13 & 14.41 \\
1s2\_4n20\_4k & 10,800    & $9.29^{}$   & 33    & 5,691.72 & 6,220.24 & 410.98    & 42.76 & 9.94  \\
1s2\_4n20\_5k & 10,800    & $6.54^{}$   & 37    & 6,169.03 & 6,572.51 & 309.37    & 38.52 & 7.12  \\
1s2\_4n25\_2k & 10,800    & N/A     & 1     & N/A      & N/A      & 10,810.15 & N/A   & N/A   \\
1s2\_4n25\_3k & 10,800    & $18.53^{}$  & 3     & 5,320.20 & 6,306.09 & 5,880.26  & 52.65 & N/A   \\
1s2\_4n25\_4k & 10,800    & $9.75^{}$   & 7     & 5,953.06 & 6,533.23 & 2,599.05  & 43.81 & N/A   \\
1s2\_4n25\_5k & 10,800    & $11.87^{}$  & 17    & 6,536.76 & 7,312.78 & 1,441.81  & 47.03 & N/A   \\
1s2\_5n5\_2k  & 22.29     & $0^{}$      & 97    & 1,712.77 & 1,712.77 & 0.99      & 56.51 & 1.64  \\
1s2\_5n5\_3k  & 90.20     & $0^{}$      & 667   & 1,897.60 & 1,897.60 & 0.64      & 48.81 & 4.87  \\
1s2\_5n5\_4k  & 4.68      & $0^{}$      & 27    & 1,923.60 & 1,923.60 & 0.51      & 34.61 & 0.18  \\
1s2\_5n5\_5k  & 80.17    & $0^{}$      & 583 & 2,228.74 & 2,228.74 & 0.41      & 41.79 & 15.95  \\
1s2\_5n10\_2k & 1,373.06  & $0^{}$      & 139   & 3,209.41 & 3,209.41 & 27.33     & 41.05 & 2.45  \\
1s2\_5n10\_3k & 559.33    & $0^{}$      & 69    & 3,471.28 & 3,471.28 & 13.56     & 39.14 & 2.88  \\
\textbf{1s2\_5n10\_4k} & \textbf{7,770.02}  & $\boldsymbol{0^{3}}$     & \textbf{1,227} & \textbf{3,855.41} & \textbf{3,855.41} & \textbf{12.66}     & \textbf{47.37} & \textbf{11.67} \\
1s2\_5n10\_5k & 886.93    & $0^{}$      & 249   & 4,083.47 & 4,083.47 & 7.65      & 41.29 & 3.86  \\
1s2\_5n15\_2k & 10,800    & $5.38^{}$   & 27    & 3,835.10 & 4,041.41 & 304.21    & 51.93 & 5.95  \\
1s2\_5n15\_3k & 10,800    & $1.91^{}$   & 95    & 4,276.14 & 4,357.70 & 125.45    & 46.50 & 3.75  \\
1s2\_5n15\_4k & 10,800    & $2.48^{3}$  & 173   & 4,473.27 & 4,584.03 & 154.69    & 41.96 & 4.97  \\
1s2\_5n15\_5k & 10,800    & $4.44^{3}$  & 171   & 4,903.99 & 5,121.54 & 111.08    & 44.83 & 8.29  \\
1s2\_5n20\_2k & 10,800    & $17.32^{}$  & 7     & 4,203.77 & 4,931.76 & 3,230.27  & 55.23 & N/A   \\
1s2\_5n20\_3k & 10,800    & $7.5^{}$    & 25    & 4,642.07 & 4,990.15 & 885.93    & 47.77 & 7.74  \\
1s2\_5n20\_4k & 10,800    & $6.18^{1}$  & 21    & 5,074.90 & 5,388.64 & 838.50    & 44.26 & 6.53  \\
1s2\_5n20\_5k & 10,800    & $3.85^{1}$  & 33    & 5,409.81 & 5,618.23 & 449.48    & 39.62 & 4.57  \\
1s2\_5n25\_2k & 10,800    & $14.91^{}$  & 3     & 4,445.92 & 5,108.63 & 8,565.08  & 50.31 & N/A   \\
1s2\_5n25\_3k & 10,800    & $8.93^{}$   & 7     & 4,948.68 & 5,390.39 & 3,163.99  & 46.14 & N/A   \\
1s2\_5n25\_4k & 10,800    & $2.61^{}$   & 15    & 5,315.99 & 5,454.63 & 1,348.20  & 36.33 & N/A   \\
1s2\_5n25\_5k & 10,800    & $5.86^{}$   & 35    & 5,714.98 & 6,049.96 & 568.92    & 40.46 & 6.43  \\* \bottomrule \\
\rlap{\textsuperscript{1} variation \#1 with parameters: $\kappa = 8$}\\
\rlap{\textsuperscript{2} variation \#2 with parameters: local depth-first search, halfpoint$ = 0.05$}\\
\rlap{\textsuperscript{3} variation \#3 with parameters: local depth-first search, $\kappa = 7$ }\\
\end{longtable}

\footnotesize  
\setlength\LTleft{-30pt}            
\setlength\LTright{-30pt} 
\begin{longtable}{@{\extracolsep{\fill}}lrrrrrrrr@{}}
\caption{Results for H3 Instances with Two Supplier and Three Satellites}
\label{tab:H3 2s3}\\
\toprule
Instance      & Time     & Gap$_f$       & \#Nodes & LB       & UB       & Time$_{root}$ & Gap$_0$  & Gap$_{20}$ \\* \midrule
\endfirsthead
\endhead
\bottomrule
\endfoot
\endlastfoot
2s3\_1n5\_2k  & 4.61     & $0^{}$      & 15      & 2,001.32 & 2,001.32 & 1.24        & 42.48  & N/A      \\
2s3\_1n5\_3k  & 14.61    & $0^{}$      & 57      & 2,119    & 2,119    & 0.60        & 37.76  & 1.75     \\
2s3\_1n5\_4k  & 12.92    & $0^{}$      & 69      & 2,232.30 & 2,232.30 & 0.74        & 33.73  & 1.26     \\
2s3\_1n5\_5k  & 46.22    & $0^{}$      & 363     & 2,477.71 & 2,477.71 & 0.41        & 40.46  & 7.75     \\
2s3\_1n10\_2k & 10,800   & $0.49^{3}$  & 403     & 3,025.89 & 3,040.69 & 39.23       & 46.61  & 7.45     \\
2s3\_1n10\_3k & 4,136.76 & $0^{}$      & 425     & 3,199.73 & 3,199.73 & 14.85       & 47.84  & 5.21     \\
2s3\_1n10\_4k & 6,296.55 & $0^{}$      & 1,075   & 3,376.60 & 3,376.60 & 9.66        & 50.26  & 14.31    \\
2s3\_1n10\_5k & 10,800   & $0.8^{3}$   & 3,311   & 3,570.39 & 3,598.89 & 6.36        & 51.24  & 7.82     \\
2s3\_1n15\_2k & 10,800   & $15.02^{3}$ & 45      & 3,471.69 & 3,993.02 & 784.45      & 52.99  & 16.76    \\
2s3\_1n15\_3k & 10,800   & $13.78^{1}$ & 59      & 3,828.31 & 4,355.79 & 228.53      & 52.38  & 16.47    \\
2s3\_1n15\_4k & 10,800   & $12.61^{3}$ & 89      & 4,107.90 & 4,626.06 & 156.58      & 50.48  & 14.88    \\
2s3\_1n15\_5k & 10,800   & $9.33^{1}$  & 143     & 4,409.15 & 4,820.51 & 97.02       & 43.71  & 11.79    \\
2s3\_1n20\_2k & 10,800   & $25.53^{}$  & 9       & 3,601.13 & 4,520.41 & 2,127.55    & 59.51  & N/A      \\
2s3\_1n20\_3k & 10,800   & $19.21^{}$  & 21      & 4,006.11 & 4,775.57 & 1,073.69    & 54.68  & 19.80    \\
2s3\_1n20\_4k & 10,800   & $12.71^{}$  & 33      & 4,333.65 & 4,884.49 & 551.02      & 49.52  & 13.49    \\
2s3\_1n20\_5k & 10,800   & $5.93^{}$   & 43      & 4,626.29 & 4,900.57 & 353.89      & 44.65  & 8.35     \\
2s3\_1n25\_2k & 10,800   & N/A         & 1       & N/A      & N/A      & 10,810.77   & N/A    & N/A      \\
2s3\_1n25\_3k & 10,800   & $22.49^{}$  & 3       & 4,314.12 & 5,284.20 & 4,712.94    & 47.28  & N/A      \\
2s3\_1n25\_4k & 10,800   & $16.85^{}$  & 15      & 4,675.25 & 5,462.81 & 1,899.50    & 42.29  & N/A      \\
2s3\_1n25\_5k & 10,800   & $12.22^{}$  & 33      & 5,005.86 & 5,617.77 & 854.59      & 36.83  & 12.71    \\
2s3\_2n5\_2k  & 363.55   & $0^{}$      & 1,437   & 1,813.62 & 1,813.62 & 1.27        & 71.12  & 13.77    \\
2s3\_2n5\_3k  & 205.34   & $0^{}$      & 1,181   & 1,877.47 & 1,877.47 & 0.83        & 75.69  & 10.71    \\
2s3\_2n5\_4k  & 255.68   & $0^{}$      & 1,839   & 2,003.37 & 2,003.37 & 0.66        & 87.48  & 13.74    \\
2s3\_2n5\_5k  & 48.12    & $0^{}$      & 361     & 2,364.94 & 2,364.94 & 0.43        & 104.25 & 25.10    \\
2s3\_2n10\_2k & 10,800   & $2.88^{3}$  & 357     & 3,695.53 & 3,802.07 & 58.78       & 51.34  & 10.88    \\
2s3\_2n10\_3k & 10,800   & $5.07^{}$   & 941     & 3,973.37 & 4,174.83 & 22.47       & 49.84  & 12.96    \\
2s3\_2n10\_4k & 10,800   & $3.75^{1}$  & 1,337   & 4,237.23 & 4,396.26 & 18.49       & 47.92  & 10.69    \\
2s3\_2n10\_5k & 10,800   & $2.57^{1}$  & 1,329   & 4,542.18 & 4,659.14 & 9.96        & 39.63  & 9        \\
2s3\_2n15\_2k & 10,800   & $11.44^{}$  & 67      & 3,545.34 & 3,950.93 & 503.42      & 61.27  & 12.51    \\
2s3\_2n15\_3k & 10,800   & $3.75^{}$   & 63      & 3,830.43 & 3,973.91 & 390.92      & 53.11  & 5.87     \\
2s3\_2n15\_4k & 10,800   & $10.85^{3}$ & 97      & 4,184.96 & 4,639.23 & 228.74      & 57.24  & 14.76    \\
2s3\_2n15\_5k & 10,800   & $3.79^{}$   & 265     & 4,462.96 & 4,632.09 & 85.35       & 45.81  & 8.05     \\
2s3\_2n20\_2k & 10,800   & $15.42^{}$  & 7       & 3,604.88 & 4,160.69 & 3,215.91    & 45.31  & N/A      \\
2s3\_2n20\_3k & 10,800   & $3.75^{}$   & 19      & 3,930.67 & 4,077.95 & 1,216.08    & 30.58  & 3.84     \\
2s3\_2n20\_4k & 10,800   & $8.31^{}$   & 23      & 4,239.05 & 4,591.49 & 654.51      & 36.39  & 8.66     \\
2s3\_2n20\_5k & 10,800   & $6.24^{1}$  & 19      & 4,567.23 & 4,852.40 & 948.16      & 31.99  & 6.30     \\
2s3\_2n25\_2k & 10,800   & $61.22^{}$  & 1       & 3,343.10 & 5,389.61 & 11,061.66   & 61.22  & N/A      \\
2s3\_2n25\_3k & 10,800   & $24.3^{}$   & 7       & 4,401.20 & 5,470.69 & 3,696.09    & 56.46  & N/A      \\
2s3\_2n25\_4k & 10,800   & $16.25^{3}$ & 9       & 4,663.69 & 5,421.55 & 2,535.76    & 46.26  & N/A      \\
2s3\_2n25\_5k & 10,800   & $7.09^{1}$  & 19      & 4,990.77 & 5,344.55 & 1,058.14    & 35.47  & 7.23     \\
2s3\_3n5\_2k  & 129.82   & $0^{}$      & 465     & 3,011.28 & 3,011.28 & 1.27        & 100.81 & 20.03    \\
2s3\_3n5\_3k  & 128.22   & $0^{}$      & 887     & 3,281.49 & 3,281.49 & 0.62        & 104.85 & 28.56    \\
2s3\_3n5\_4k  & 32.88    & $0^{}$      & 209     & 3,430.70 & 3,430.70 & 0.53        & 88.15  & 22.04    \\
2s3\_3n5\_5k  & 5.90     & $0^{}$      & 69      & 3,573.10 & 3,573.10 & 0.44        & 76.39  & 19.50    \\
2s3\_3n10\_2k & 1,902.44 & $0^{}$      & 91      & 2,706.42 & 2,706.42 & 61.95       & 64.90  & 2.20     \\
2s3\_3n10\_3k & 339.32   & $0^{}$      & 27      & 2,792.60 & 2,792.60 & 40.16       & 51.92  & 1.34     \\
\textbf{2s3\_3n10\_4k} & \textbf{6,503.28} & $\boldsymbol{0^{3}}$     & \textbf{837}     & \textbf{3,070.92} & \textbf{3,070.92} & \textbf{31.21}       & \textbf{51.68}  & \textbf{6.06}     \\
2s3\_3n10\_5k & 5,422.97 & $0^{}$      & 1,609   & 3,247    & 3,247    & 12.14       & 45.66  & 5.21     \\
2s3\_3n15\_2k & 10,800   & $13.72^{}$  & 61      & 3,919.09 & 4,456.83 & 520.35      & 51.21  & 16.37    \\
2s3\_3n15\_3k & 10,800   & $11.92^{1}$ & 51      & 4,279.70 & 4,789.84 & 592.84      & 53.29  & 13.42    \\
2s3\_3n15\_4k & 10,800   & $6.87^{3}$  & 125     & 4,601.45 & 4,917.43 & 133.66      & 49.70  & 12.02    \\
2s3\_3n15\_5k & 10,800   & $6.19^{}$   & 413     & 4,907.58 & 5,211.57 & 73.45       & 47.90  & 11.68    \\
2s3\_3n20\_2k & 10,800   & $21.02^{3}$ & 5       & 3,874.73 & 4,689.01 & 2,982.58    & 53.63  & N/A      \\
2s3\_3n20\_3k & 10,800   & $17.68^{}$  & 27      & 4,109.85 & 4,836.30 & 675.72      & 50.85  & 18.17    \\
2s3\_3n20\_4k & 10,800   & $3.87^{3}$  & 21      & 4,458.92 & 4,631.50 & 696.87      & 34.09  & 4.13     \\
2s3\_3n20\_5k & 10,800   & $2.75^{3}$  & 37      & 4,710.66 & 4,840.12 & 446.44      & 30.46  & 4.09     \\
2s3\_3n25\_2k & 10,800   & $47.29^{}$  & 1       & 3,244.50 & 4,778.68 & 8,779.68    & 47.29  & N/A      \\
2s3\_3n25\_3k & 10,800   & $16.91^{}$  & 9       & 4,286.28 & 5,011.15 & 2,444.44    & 48     & N/A      \\
2s3\_3n25\_4k & 10,800   & $12.17^{}$  & 19      & 4,579.88 & 5,137.33 & 1,240.57    & 41.19  & 12.28    \\
2s3\_3n25\_5k & 10,800   & $8.19^{}$   & 27      & 4,803.79 & 5,197.39 & 709.51      & 34.90  & 8.58     \\
2s3\_4n5\_2k  & 64.61    & $0^{}$      & 131     & 2,247.55 & 2,247.55 & 1.55        & 49.41  & 10.07    \\
2s3\_4n5\_3k  & 92.91    & $0^{}$      & 453     & 2,563.67 & 2,563.67 & 0.77        & 53.01  & 13.65    \\
2s3\_4n5\_4k  & 24.51    & $0^{}$      & 167     & 2,676.43 & 2,676.43 & 0.55        & 40.82  & 6.75     \\
2s3\_4n5\_5k  & 965.71    & $0^{}$      & 6653     & 3,327.69 & 3,327.69 & 0.45        & 56.61  & 27.18    \\
2s3\_4n10\_2k & 10,800   & $3.6^{}$    & 643     & 3,321.15 & 3,440.56 & 38.65       & 53.75  & 9.68     \\
2s3\_4n10\_3k & 10,800   & $6.93^{}$   & 1,355   & 3,633.12 & 3,884.92 & 18.72       & 63.44  & 16.45    \\
2s3\_4n10\_4k & 10,800   & $6.23^{}$   & 1,847   & 3,966.92 & 4,214.22 & 11.33       & 62.98  & 18.62    \\
2s3\_4n10\_5k & 10,800   & $5.96^{}$   & 3,047   & 4,283.86 & 4,539.38 & 6.94        & 64     & 21.12    \\
2s3\_4n15\_2k & 10,800   & $18.03^{}$  & 25      & 3,586.75 & 4,233.31 & 1,068.46    & 59.08  & 18.33    \\
2s3\_4n15\_3k & 10,800   & $8.63^{}$   & 45      & 4,151.37 & 4,509.53 & 558.23      & 48.60  & 10.83    \\
2s3\_4n15\_4k & 10,800   & $5^{}$      & 49      & 4,468.48 & 4,692.03 & 408.68      & 38.55  & 6.64     \\
2s3\_4n15\_5k & 10,800   & $2.16^{}$   & 69      & 4,837.24 & 4,941.65 & 254.61      & 29.60  & 3.62     \\
2s3\_4n20\_2k & 10,800   & $30.27^{}$  & 11      & 3,960.31 & 5,158.97 & 1,896.29    & 62.58  & N/A      \\
2s3\_4n20\_3k & 10,800   & $18.58^{3}$ & 7       & 4,261.53 & 5,053.28 & 1,591.58    & 45.69  & N/A      \\
2s3\_4n20\_4k & 10,800   & $14.06^{}$  & 37      & 4,739.12 & 5,405.30 & 492.08      & 42.96  & 16.42    \\
2s3\_4n20\_5k & 10,800   & $8.18^{1}$  & 31      & 5,150.20 & 5,571.28 & 603.81      & 33.25  & 9        \\
2s3\_4n25\_2k & 10,800   & N/A         & 1       & N/A      & N/A      & 10,810.29   & N/A    & N/A      \\
2s3\_4n25\_3k & 10,800   & $29.09^{}$  & 5       & 4,534.46 & 5,853.40 & 3,845.42    & 62.57  & N/A      \\
2s3\_4n25\_4k & 10,800   & $17.27^{}$  & 13      & 4,913.74 & 5,762.48 & 2,049.71    & 46.85  & N/A      \\
2s3\_4n25\_5k & 10,800   & $10.04^{}$  & 19      & 5,224.12 & 5,748.56 & 1,296.60    & 36.68  & 10.06    \\
2s3\_5n5\_2k  & 157.24   & $0^{}$      & 693     & 1,685.23 & 1,685.23 & 1.39        & 72.47  & 12.78    \\
2s3\_5n5\_3k  & 35.53    & $0^{}$      & 223     & 1,780.87 & 1,780.87 & 0.82        & 61.01  & 6.06     \\
2s3\_5n5\_4k  & 24.35    & $0^{}$      & 223     & 1,879.77 & 1,879.77 & 0.54        & 53.31  & 3.58     \\
2s3\_5n5\_5k  & 199.95   & $0^{}$      & 1,913   & 2,036.02 & 2,036.02 & 0.37        & 58.04  & 11.99    \\
2s3\_5n10\_2k & 10,800   & $1.23^{}$   & 1,043   & 3,081.67 & 3,119.51 & 20.17       & 46.13  & 9.35     \\
2s3\_5n10\_3k & 3,069.23 & $0^{}$      & 437     & 3,246.90 & 3,246.90 & 23.06       & 42.08  & 6.80     \\
2s3\_5n10\_4k & 10,800   & $3.31^{}$   & 2,199   & 3,456.33 & 3,570.82 & 18.09       & 50.12  & 13.34    \\
2s3\_5n10\_5k & 2,456.68 & $0^{}$      & 663     & 3,673.39 & 3,673.39 & 6.78        & 48.22  & 16.01    \\
2s3\_5n15\_2k & 10,800   & $9.34^{3}$  & 29      & 3,400.87 & 3,718.53 & 623.70      & 46.20  & 10.10    \\
2s3\_5n15\_3k & 10,800   & $9.98^{1}$  & 31      & 3,752.60 & 4,127.05 & 331.19      & 48.85  & 11.80    \\
2s3\_5n15\_4k & 10,800   & $13.14^{1}$ & 153     & 3,892.41 & 4,403.72 & 150.96      & 55.53  & 17.18    \\
2s3\_5n15\_5k & 10,800   & $9.39^{3}$  & 189     & 4,355.31 & 4,764.33 & 78.03       & 54.10  & 17.22    \\
2s3\_5n20\_2k & 10,800   & $19.98^{}$  & 5       & 3,726.31 & 4,470.75 & 4,515.41    & 51.78  & N/A      \\
2s3\_5n20\_3k & 10,800   & $13.76^{}$  & 25      & 4,067.21 & 4,626.79 & 918.42      & 44.94  & 15.14    \\
2s3\_5n20\_4k & 10,800   & $13.72^{}$  & 35      & 4,419.30 & 5,025.79 & 665.38      & 45.08  & 14.98    \\
2s3\_5n20\_5k & 10,800   & $6.49^{1}$  & 29      & 4,670.17 & 4,973.15 & 559.87      & 34.22  & 8.05     \\
2s3\_5n25\_2k & 10,800   & $27.88^{}$  & 3       & 3,930.41 & 5,026.24 & 5,563.11    & 60.61  & N/A      \\
2s3\_5n25\_3k & 10,800   & $23.67^{1}$ & 3       & 4,127.03 & 5,103.80 & 4,846.97    & 54.23  & N/A      \\
2s3\_5n25\_4k & 10,800   & $6.63^{}$   & 23      & 4,686.88 & 4,997.70 & 837.07      & 45.11  & 10.93    \\
2s3\_5n25\_5k & 10,800   & $8.79^{}$   & 35      & 4,999.72 & 5,439.04 & 600.31      & 47.58  & 16.09    \\* \bottomrule\\
\rlap{\textsuperscript{1} variation \#1 with parameters: $\kappa = 8$}\\
\rlap{\textsuperscript{2} variation \#2 with parameters: local depth-first search, halfpoint$ = 0.05$}\\
\rlap{\textsuperscript{3} variation \#3 with parameters: local depth-first search, $\kappa = 7$ }\\
\end{longtable}
\footnotesize  
\setlength\LTleft{-30pt}            
\setlength\LTright{-30pt} 
\begin{longtable}{@{\extracolsep{\fill}}lrrrrrrrr@{}}
\caption{Results for L3 Instances with One Supplier and Two Satellites}
\label{tab:L3 1s2}\\
\toprule
Instance &
  \multicolumn{1}{l}{Time} &
  \multicolumn{1}{l}{Gap$_f$} &
  \multicolumn{1}{l}{\#Nodes} &
  \multicolumn{1}{l}{LB} &
  \multicolumn{1}{l}{UB} &
  \multicolumn{1}{l}{Time$_{root}$} &
  \multicolumn{1}{l}{Gap$_0$} &
  \multicolumn{1}{l}{Gap$_{20}$} \\* \midrule
\endfirsthead
\endhead
\bottomrule
\endfoot
\endlastfoot
1s2\_1n5\_2k  & 23.95    & $0^{}$      & 87    & 2,684.71 & 2,684.71 & 1.16      & 56.84 & 3.61  \\
1s2\_1n5\_3k  & 65.12    & $0^{}$      & 355   & 3,011.36 & 3,011.36 & 0.87      & 47.44 & 5.17  \\
1s2\_1n5\_4k  & 43.23    & $0^{}$      & 245   & 3,218.80 & 3,218.80 & 0.73      & 34.78 & 3.50  \\
1s2\_1n5\_5k  & 79.94    & $0^{}$      & 749   & 3,452.14 & 3,452.14 & 0.40      & 32.58 & 3.13  \\
1s2\_1n10\_2k & 8,828.05 & $0^{}$      & 599   & 3,044    & 3,044    & 52.47     & 57.65 & 7.30  \\
1s2\_1n10\_3k & 6,960.42 & $0^{}$      & 945   & 3,284.53 & 3,284.53 & 17.67     & 49.63 & 6.15  \\
\textbf{1s2\_1n10\_4k} & \textbf{9,647.56} & $\boldsymbol{0^{2}}$     & \textbf{1,739} & \textbf{3,545.64} & \textbf{3,545.64} & \textbf{10.03}     & \textbf{50.11} & \textbf{2.99}  \\
1s2\_1n10\_5k & 10,800   & $3.7^{1}$   & 2,559 & 3,768.75 & 3,908.15 & 7.84      & 46.48 & 6.45  \\
1s2\_1n15\_2k & 10,800   & $4.99^{}$   & 69    & 3,675.24 & 3,858.47 & 397.37    & 55.44 & 6.43  \\
1s2\_1n15\_3k & 10,800   & $5.47^{1}$  & 95    & 4,181.14 & 4,409.64 & 259.73    & 47.53 & 7.53  \\
1s2\_1n15\_4k & 10,800   & $5^{1}$     & 131   & 4,564.85 & 4,793.08 & 305.20    & 41.26 & 6.73  \\
1s2\_1n15\_5k & 10,800   & $8.46^{1}$  & 123   & 4,998.14 & 5,420.75 & 122.81    & 43.57 & 10.51 \\
1s2\_1n20\_2k & 10,800   & $16.69^{3}$ & 3     & 3,434.56 & 4,007.96 & 6,515.97  & 65.58 & N/A   \\
1s2\_1n20\_3k & 10,800   & $5.1^{}$    & 25    & 3,872.30 & 4,069.77 & 1,102.05  & 51.52 & 5.41  \\
1s2\_1n20\_4k & 10,800   & $7.15^{1}$  & 19    & 4,136.95 & 4,432.89 & 1,201.17  & 48.95 & 7.21  \\
1s2\_1n20\_5k & 10,800   & $0.91^{}$   & 85    & 4,510.96 & 4,551.97 & 241.55    & 41.63 & 2.50  \\
1s2\_1n25\_2k & 10,800   & $19.96^{}$  & 3     & 3,858.27 & 4,628.24 & 7,565.21  & 64.79 & N/A   \\
1s2\_1n25\_3k & 10,800   & $19.8^{1}$  & 3     & 4,245.69 & 5,086.24 & 6,461.65  & 58.93 & N/A   \\
1s2\_1n25\_4k & 10,800   & $13.74^{}$  & 15    & 4,595.73 & 5,227.23 & 1,488.47  & 54.70 & N/A   \\
1s2\_1n25\_5k & 10,800   & $3.79^{}$   & 45    & 5,049.49 & 5,241.10 & 582.75    & 40.14 & 4.23  \\
1s2\_2n5\_2k  & 70.53    & $0^{}$      & 373   & 2,225.84 & 2,225.84 & 1.11      & 77.98 & 7.97  \\
1s2\_2n5\_3k  & 20.65    & $0^{}$      & 165   & 2,299.56 & 2,299.56 & 0.61      & 67.25 & 3.52  \\
1s2\_2n5\_4k  & 25.90    & $0^{}$      & 171   & 2,480.16 & 2,480.16 & 0.57      & 72.28 & 3.42  \\
1s2\_2n5\_5k  & 73.59    & $0^{}$      & 511   & 2,708.57 & 2,708.57 & 0.44      & 73.67 & 4.41  \\
1s2\_2n10\_2k & 10,800   & $3.31^{}$   & 643   & 3,917.01 & 4,046.60 & 18.34     & 74.39 & 11.18 \\
1s2\_2n10\_3k & 4,693.68 & $0^{}$      & 365   & 4,294.46 & 4,294.46 & 17.50     & 60.41 & 8.07  \\
1s2\_2n10\_4k & 10,800   & $4.58^{}$   & 2,471 & 4,691.96 & 4,906.87 & 13.80     & 66.31 & 11.33 \\
1s2\_2n10\_5k & 5,048.29 & $0^{}$      & 1,353 & 4,997.90 & 4,997.90 & 11.03     & 49.37 & 6.35  \\
1s2\_2n15\_2k & 10,800   & $7.4^{3}$   & 23    & 3,813.91 & 4,096.28 & 814.61    & 62.58 & 7.99  \\
1s2\_2n15\_3k & 10,800   & $6.31^{1}$  & 45    & 4,181.02 & 4,444.76 & 487.76    & 53.55 & 7.60  \\
1s2\_2n15\_4k & 10,800   & $6.17^{3}$  & 83    & 4,793.53 & 5,089.53 & 216.71    & 49.90 & 8.50  \\
1s2\_2n15\_5k & 10,800   & $3.66^{}$   & 221   & 5,192.10 & 5,382.14 & 76.21     & 42.24 & 5.88  \\
1s2\_2n20\_2k & 10,800   & $16.03^{}$  & 11    & 3,968.86 & 4,605.20 & 2,636.03  & 72.19 & N/A   \\
1s2\_2n20\_3k & 10,800   & $9.27^{}$   & 29    & 4,403.94 & 4,811.99 & 780.71    & 59.47 & 9.53  \\
1s2\_2n20\_4k & 10,800   & $4.42^{1}$  & 19    & 4,784.34 & 4,995.98 & 901.73    & 49.44 & 4.44  \\
1s2\_2n20\_5k & 10,800   & $4.41^{3}$  & 25    & 5,203.90 & 5,433.15 & 476.30    & 48.06 & 4.73  \\
1s2\_2n25\_2k & 10,800   & N/A         & 1     & N/A      & N/A      & 10,810.36 & N/A   & N/A   \\
1s2\_2n25\_3k & 10,800   & $19.85^{1}$ & 3     & 4,276.60 & 5,125.57 & 6,262.92  & 61.61 & N/A   \\
1s2\_2n25\_4k & 10,800   & $13.57^{}$  & 13    & 4,912.72 & 5,579.59 & 1,446.51  & 64.79 & N/A   \\
1s2\_2n25\_5k & 10,800   & $5.39^{}$   & 21    & 5,387.19 & 5,677.72 & 839.23    & 50.45 & 5.53  \\
1s2\_3n5\_2k  & 50.19    & $0^{}$      & 231   & 3,278.95 & 3,278.95 & 0.78      & 71.87 & 5.14  \\
1s2\_3n5\_3k  & 65.86    & $0^{}$      & 399   & 3,655.75 & 3,655.75 & 0.55      & 83.26 & 9.71  \\
1s2\_3n5\_4k  & 238.21   & $0^{}$      & 1,531 & 4,101.64 & 4,101.64 & 0.59      & 76.08 & 8.20  \\
1s2\_3n5\_5k  & 12.90    & $0^{}$      & 127   & 4,141.28 & 4,141.28 & 0.36      & 64.60 & 1.81  \\
1s2\_3n10\_2k & 10,800   & $4.93^{}$   & 181   & 3,302.71 & 3,465.57 & 119.57    & 62.94 & 8.20  \\
1s2\_3n10\_3k & 10,800   & $2.23^{3}$  & 553   & 3,749.96 & 3,833.46 & 23.46     & 54.02 & 6.65  \\
1s2\_3n10\_4k & 10,800   & $7.31^{1}$  & 481   & 4,228.54 & 4,537.52 & 26.89     & 54.11 & 8.97  \\
1s2\_3n10\_5k & 10,800   & $4.67^{1}$  & 903   & 4,633.44 & 4,850.03 & 13.49     & 49.88 & 6.75  \\
1s2\_3n15\_2k & 10,800   & $6.5^{1}$   & 23    & 4,076.08 & 4,341    & 935.07    & 58.33 & 6.82  \\
1s2\_3n15\_3k & 10,800   & $2.34^{}$   & 133   & 4,522.98 & 4,628.73 & 214.71    & 49.72 & 3.88  \\
1s2\_3n15\_4k & 10,800   & $2.57^{}$   & 217   & 4,793.54 & 4,916.73 & 104.64    & 48.57 & 5.28  \\
1s2\_3n15\_5k & 3,697.07 & $0^{}$      & 109   & 5,126.95 & 5,126.95 & 70.55     & 37.13 & 1.13  \\
1s2\_3n20\_2k & 10,800   & $22.23^{1}$ & 3     & 3,777.01 & 4,616.67 & 5,923.74  & 71.84 & N/A   \\
1s2\_3n20\_3k & 10,800   & $9.82^{}$   & 27    & 4,201.56 & 4,614.29 & 1,048.46  & 50.41 & 10.07 \\
1s2\_3n20\_4k & 10,800   & $6.21^{}$   & 49    & 4,689.40 & 4,980.59 & 546.25    & 42.89 & 7.35  \\
1s2\_3n20\_5k & 10,800   & $5.33^{1}$  & 23    & 5,101.31 & 5,373.19 & 777.79    & 36.58 & 5.56  \\
1s2\_3n25\_2k & 10,800   & N/A         & 1     & N/A      & N/A      & 10,801.88 & N/A   & N/A   \\
1s2\_3n25\_3k & 10,800   & $12.77^{}$  & 7     & 4,271.40 & 4,816.81 & 3,417.78  & 57.33 & N/A   \\
1s2\_3n25\_4k & 10,800   & $9.93^{}$   & 11    & 4,672.46 & 5,136.39 & 1,276     & 56.63 & N/A   \\
1s2\_3n25\_5k & 10,800   & $5.25^{}$   & 21    & 4,944.52 & 5,203.89 & 971.59    & 49.21 & 5.35  \\
1s2\_4n5\_2k  & 40.08    & $0^{}$      & 107   & 2,575.53 & 2,575.53 & 1.34      & 56.09 & 8.24  \\
1s2\_4n5\_3k  & 39.11    & $0^{}$      & 189   & 2,957.49 & 2,957.49 & 0.64      & 48.49 & 6.32  \\
1s2\_4n5\_4k  & 207.04   & $0^{}$      & 1,211 & 3,317.46 & 3,317.46 & 0.79      & 36.68 & 4.97  \\
1s2\_4n5\_5k  & 17.20     & $0^{}$      & 133    & 3,726.02 & 3,726.02 & 0.30      & 30.49 & 2.98  \\
1s2\_4n10\_2k & 3,108.44 & $0^{}$      & 181   & 3,785.20 & 3,785.20 & 37.15     & 49.43 & 4.62  \\
\textbf{1s2\_4n10\_3k} & \textbf{8,472.44} & $\boldsymbol{0^{3}}$     & \textbf{895}   & \textbf{4,350.14} & \textbf{4,350.14} & \textbf{24.42}     & \textbf{49.08} & \textbf{9.97}  \\
1s2\_4n10\_4k & 2,403.62 & $0^{}$      & 461   & 4,697.63 & 4,697.63 & 11.81     & 45.98 & 5.62  \\
1s2\_4n10\_5k & 10,800   & $5.74^{}$   & 3,351 & 5,159.91 & 5,456.08 & 6.49      & 49.79 & 12.63 \\
1s2\_4n15\_2k & 10,800   & $9.88^{}$   & 25    & 3,388.76 & 3,723.47 & 933.96    & 63.46 & 10.25 \\
1s2\_4n15\_3k & 10,800   & $3.51^{}$   & 41    & 3,676.51 & 3,805.55 & 376.37    & 45.63 & 4.96  \\
1s2\_4n15\_4k & 10,800   & $2.54^{1}$  & 35    & 3,953.66 & 4,053.91 & 437.83    & 37.45 & 3.33  \\
1s2\_4n15\_5k & 10,800   & $2.96^{}$   & 103   & 4,241.38 & 4,366.88 & 154.72    & 32.48 & 3.91  \\
1s2\_4n20\_2k & 10,800   & $18.2^{1}$  & 9     & 4,292.96 & 5,074.49 & 2,981.84  & 62.43 & N/A   \\
1s2\_4n20\_3k & 10,800   & $15.43^{}$  & 27    & 4,768.58 & 5,504.58 & 695.47    & 56.57 & 16.04 \\
1s2\_4n20\_4k & 10,800   & $8.51^{1}$  & 17    & 5,211.76 & 5,655.09 & 1,295.73  & 43.42 & N/A   \\
1s2\_4n20\_5k & 10,800   & $5.38^{1}$  & 33    & 5,699.97 & 6,006.71 & 606.04    & 37.75 & 6.13  \\
1s2\_4n25\_2k & 10,800   & N/A         & 1     & N/A      & N/A      & 10,810.32 & N/A   & N/A   \\
1s2\_4n25\_3k & 10,800   & $21.29^{}$  & 5     & 4,783.63 & 5,801.83 & 5,353.92  & 63.80 & N/A   \\
1s2\_4n25\_4k & 10,800   & $14.7^{}$   & 9     & 5,210.68 & 5,976.54 & 2,562.62  & 51.08 & N/A   \\
1s2\_4n25\_5k & 10,800   & $10.58^{}$  & 21    & 5,941.68 & 6,570.19 & 1,044.59  & 49.41 & 10.67 \\
1s2\_5n5\_2k  & 8.21     & $0^{}$      & 33    & 1,610.91 & 1,610.91 & 0.99      & 57.88 & 0.93  \\
1s2\_5n5\_3k  & 189.67   & $0^{}$      & 1,789 & 1,782.25 & 1,782.25 & 0.54      & 48.04 & 5.96  \\
1s2\_5n5\_4k  & 61.82    & $0^{}$      & 475   & 1,852.80 & 1,852.80 & 0.63      & 35.81 & 1.11  \\
1s2\_5n5\_5k  & 540.40   & $0^{}$      & 4,375 & 2,095.40 & 2,095.40 & 0.48      & 44.66 & 9.41  \\
1s2\_5n10\_2k & 1,119.83 & $0^{}$      & 167   & 2,848.08 & 2,848.08 & 26.44     & 44.35 & 2.35  \\
1s2\_5n10\_3k & 802.09   & $0^{}$      & 125   & 3,123.52 & 3,123.52 & 15.77     & 43.92 & 2.32  \\
\textbf{1s2\_5n10\_4k} & \textbf{9,394.96} & $\boldsymbol{0^{1}}$     & \textbf{931}   & \textbf{3,494.08} & \textbf{3,494.08} & \textbf{16.91}     & \textbf{51.85} & \textbf{8.49}  \\
1s2\_5n10\_5k & 550.34   & $0^{}$      & 181   & 3,735.42 & 3,735.42 & 7.03      & 42.89 & 4.23  \\
1s2\_5n15\_2k & 10,800   & $7.53^{}$   & 53    & 3,437.19 & 3,695.90 & 332.90    & 60.65 & 8.61  \\
1s2\_5n15\_3k & 10,800   & $1.22^{}$   & 159   & 3,888.06 & 3,935.50 & 140.66    & 50.08 & 3.23  \\
1s2\_5n15\_4k & 10,800   & $0.63^{}$   & 323   & 4,072.26 & 4,098.06 & 62.16     & 44.07 & 5.21  \\
1s2\_5n15\_5k & 10,800   & $3.45^{3}$  & 227   & 4,509.66 & 4,665.03 & 55.99     & 46.10 & 8.99  \\
1s2\_5n20\_2k & 10,800   & $17.17^{}$  & 11    & 3,625.51 & 4,247.96 & 3,163.87  & 63.61 & N/A   \\
1s2\_5n20\_3k & 10,800   & $9.73^{}$   & 33    & 4,068.06 & 4,463.70 & 836.97    & 59.14 & 10.10 \\
1s2\_5n20\_4k & 10,800   & $6.21^{3}$  & 25    & 4,470.60 & 4,748.28 & 813.28    & 48.74 & 6.61  \\
1s2\_5n20\_5k & 10,800   & $5^{3}$     & 37    & 4,820.07 & 5,061.15 & 285.11    & 44.68 & 5.90  \\
1s2\_5n25\_2k & 10,800   & $64.27^{}$  & 1     & 2,629    & 4,318.71 & 10,225.16 & 64.27 & N/A   \\
1s2\_5n25\_3k & 10,800   & $6.73^{}$   & 11    & 4,182.78 & 4,464.21 & 2,156.88  & 49.75 & N/A   \\
1s2\_5n25\_4k & 4,575.74 & $0^{}$      & 11    & 4,524.71 & 4,524.71 & 942.34    & 37.94 & N/A   \\
1s2\_5n25\_5k & 10,800   & $5.79^{3}$  & 25    & 4,943.96 & 5,230.23 & 909.45    & 44.29 & 6.14  \\* \bottomrule\\
\rlap{\textsuperscript{1} variation \#1 with parameters: $\kappa = 8$}\\
\rlap{\textsuperscript{2} variation \#2 with parameters: local depth-first search, halfpoint$ = 0.05$}\\
\rlap{\textsuperscript{3} variation \#3 with parameters: local depth-first search, $\kappa = 7$ }\\
\end{longtable}

\footnotesize  
\setlength\LTleft{-30pt}            
\setlength\LTright{-30pt} 
\begin{longtable}{@{\extracolsep{\fill}}lrrrrrrrr@{}}
\caption{Results for L3 Instances with Two Suppliers and Three Satellites}
\label{tab:L3 2s3}\\
\toprule
Instance &
  \multicolumn{1}{l}{Time} &
  \multicolumn{1}{l}{Gap$_f$} &
  \multicolumn{1}{l}{\#Nodes} &
  \multicolumn{1}{l}{LB} &
  \multicolumn{1}{l}{UB} &
  \multicolumn{1}{l}{Time$_{root}$} &
  \multicolumn{1}{l}{Gap$_0$} &
  \multicolumn{1}{l}{Gap$_{20}$} \\* \midrule
\endfirsthead
\endhead
\bottomrule
\endfoot
\endlastfoot
2s3\_1n5\_2k  & 4.39      & $0^{}$      & 11     & 1,879.50 & 1,879.50 & 1.37      & 38.94  & N/A   \\
2s3\_1n5\_3k  & 14.01     & $0^{}$      & 53     & 2,001.02 & 2,001.02 & 0.75      & 35.92  & 3.06  \\
2s3\_1n5\_4k  & 6.84      & $0^{}$      & 43     & 2,115.85 & 2,115.85 & 0.63      & 31.13  & 1.27  \\
2s3\_1n5\_5k  & 78.28     & $0^{}$      & 593    & 2,366.88 & 2,366.88 & 0.38      & 38.12  & 8.32  \\
2s3\_1n10\_2k & 10,800    & $3.36^{3}$  & 555    & 2,709.78 & 2,800.75 & 45.52     & 54.85  & 11.06 \\
2s3\_1n10\_3k & 2,461.06  & $0^{}$      & 279    & 2,896.02 & 2,896.02 & 19.33     & 54.35  & 10.82 \\
2s3\_1n10\_4k & 587.67    & $0^{}$      & 113    & 3,010.55 & 3,010.55 & 9.24      & 52.23  & 21.99 \\
2s3\_1n10\_5k & 10,800    & $3.85^{1}$  & 2,637  & 3,302.38 & 3,429.63 & 5.17      & 59.53  & 16.30 \\
2s3\_1n15\_2k & 10,800    & $8.58^{1}$  & 19     & 3,064.06 & 3,327.05 & 1,242.78  & 47.96  & 8.60  \\
2s3\_1n15\_3k & 10,800    & $12.79^{3}$ & 47     & 3,459.82 & 3,902.43 & 224.25    & 54.26  & 14.62 \\
2s3\_1n15\_4k & 10,800    & $11.54^{3}$ & 91     & 3,771.38 & 4,206.67 & 138.01    & 53.01  & 14.93 \\
2s3\_1n15\_5k & 10,800    & $10.73^{1}$ & 85     & 4,047.21 & 4,481.30 & 99.28     & 47.52  & 14.05 \\
2s3\_1n20\_2k & 10,800    & $23.49^{}$  & 9      & 3,065.53 & 3,785.50 & 3,690.57  & 61.27  & N/A   \\
2s3\_1n20\_3k & 10,800    & $18.19^{}$  & 21     & 3,525.40 & 4,166.81 & 1,146.70  & 60.05  & 18.50 \\
2s3\_1n20\_4k & 10,800    & $10.26^{}$  & 31     & 3,888.10 & 4,287.17 & 505.19    & 53.65  & 13.41 \\
2s3\_1n20\_5k & 10,800    & $10.27^{3}$ & 35     & 4,219.07 & 4,652.55 & 503.24    & 59.75  & 15.05 \\
2s3\_1n25\_2k & 10,800    & N/A         & 1      & N/A      & N/A      & 10,810.01 & N/A    & N/A   \\
2s3\_1n25\_3k & 10,800    & $28.59^{}$  & 5      & 3,721.68 & 4,785.56 & 3,705.21  & 57.53  & N/A   \\
2s3\_1n25\_4k & 10,800    & $26.27^{}$  & 17     & 4,091.91 & 5,166.87 & 2,174.16  & 56.99  & N/A   \\
2s3\_1n25\_5k & 10,800    & $16.08^{}$  & 35     & 4,477.76 & 5,197.72 & 1,145.69  & 45.94  & 16.79 \\
2s3\_2n5\_2k  & 224.65    & $0^{}$      & 1,177  & 1,760.91 & 1,760.91 & 0.98      & 71.93  & 11.49 \\
2s3\_2n5\_3k  & 135.56    & $0^{}$      & 973    & 1,824.03 & 1,824.03 & 0.77      & 75.52  & 8.02  \\
2s3\_2n5\_4k  & 863.65    & $0^{}$      & 5,749  & 1,993.91 & 1,993.91 & 0.54      & 91.82  & 14.11 \\
2s3\_2n5\_5k  & 148.15    & $0^{}$      & 1,107  & 2,332.61 & 2,332.61 & 0.45      & 106.59 & 27.18 \\
\textbf{2s3\_2n10\_2k} & \textbf{7,579.55}  & $\boldsymbol{0^{1}}$     & \textbf{229}    & \textbf{3,407.36} & \textbf{3,407.36} & \textbf{52.14}     & \textbf{49.32}  & \textbf{7.30}  \\
2s3\_2n10\_3k & 7,247.46  & $0^{}$      & 865    & 3,721.08 & 3,721.08 & 14.54     & 44.82  & 9.04  \\
2s3\_2n10\_4k & 10,800    & $8.55^{1}$  & 1,809  & 3,958.46 & 4,296.75 & 14.55     & 56.37  & 19.48 \\
2s3\_2n10\_5k & 10,800    & $1.73^{1}$  & 1,625  & 4,270.69 & 4,344.44 & 6.45      & 37.82  & 7.83  \\
2s3\_2n15\_2k & 10,800    & $11.61^{1}$ & 21     & 3,231.75 & 3,607.06 & 1,962.59  & 66.91  & 14.12 \\
2s3\_2n15\_3k & 10,800    & $16.95^{}$  & 67     & 3,493.47 & 4,085.74 & 297.60    & 81.52  & 20.99 \\
2s3\_2n15\_4k & 10,800    & $12.45^{}$  & 143    & 3,889.94 & 4,374.19 & 170.58    & 67.37  & 18.45 \\
2s3\_2n15\_5k & 10,800    & $11.94^{1}$ & 185    & 4,181.31 & 4,680.73 & 111.77    & 64.90  & 19.70 \\
2s3\_2n20\_2k & 10,800    & $22.45^{}$  & 5      & 3,096.44 & 3,791.61 & 3,404.59  & 58.65  & N/A   \\
2s3\_2n20\_3k & 10,800    & $18.44^{}$  & 19     & 3,455.67 & 4,092.83 & 1,157.40  & 53.06  & 18.49 \\
2s3\_2n20\_4k & 10,800    & $10.02^{}$  & 33     & 3,801.63 & 4,182.53 & 403.43    & 43.01  & 10.66 \\
2s3\_2n20\_5k & 10,800    & $13.58^{3}$ & 31     & 4,161.01 & 4,726.12 & 443.42    & 45.18  & 14.02 \\
2s3\_2n25\_2k & 10,800    & $63.54^{}$  & 1      & 2,684.59 & 4,390.27 & 10,280.58 & 63.54  & N/A   \\
2s3\_2n25\_3k & 10,800    & $27.53^{}$  & 11     & 3,713.56 & 4,736.04 & 2,814.18  & 66.28  & N/A   \\
2s3\_2n25\_4k & 10,800    & $11.6^{}$   & 21     & 3,994.76 & 4,458.23 & 1,467.72  & 46.11  & 11.81 \\
2s3\_2n25\_5k & 10,800    & $6.31^{}$   & 39     & 4,361.90 & 4,637.05 & 802.44    & 39.64  & 9.67  \\
2s3\_3n5\_2k  & 35.95     & $0^{}$      & 155    & 2,793.06 & 2,793.06 & 1.06      & 102.87 & 22.71 \\
2s3\_3n5\_3k  & 92.66     & $0^{}$      & 555    & 3,088.87 & 3,088.87 & 0.67      & 109.31 & 27.70 \\
2s3\_3n5\_4k  & 13.35     & $0^{}$      & 91     & 3,234.17 & 3,234.17 & 0.60      & 85.61  & 16.91 \\
2s3\_3n5\_5k  & 12.32     & $0^{}$      & 111    & 3,379.81 & 3,379.81 & 0.41      & 73.21  & 16.49 \\
2s3\_3n10\_2k & 941       & $0^{}$      & 61     & 2,434.38 & 2,434.38 & 54        & 63.90  & 1.71  \\
2s3\_3n10\_3k & 278.78    & $0^{}$      & 41     & 2,521.19 & 2,521.19 & 20.60     & 49.99  & 1.73  \\
2s3\_3n10\_4k & 10,800    & $0.85^{}$   & 2,379  & 2,779.42 & 2,802.93 & 18.68     & 50.33  & 9.48  \\
2s3\_3n10\_5k & 10,479.90 & $0^{}$      & 2,509  & 2,972.77 & 2,972.77 & 12.83     & 43.21  & 4.99  \\
2s3\_3n15\_2k & 10,800    & $9.78^{3}$  & 33     & 3,474.34 & 3,814.29 & 1,697.92  & 46.77  & 10.44 \\
2s3\_3n15\_3k & 10,800    & $6.33^{3}$  & 75     & 3,909.93 & 4,157.29 & 607.33    & 52.06  & 8.19  \\
2s3\_3n15\_4k & 10,800    & $8.76^{1}$  & 175    & 4,236.42 & 4,607.55 & 122.57    & 57.34  & 14.54 \\
2s3\_3n15\_5k & 10,800    & $5.5^{1}$   & 129    & 4,493    & 4,740.06 & 72.03     & 44.86  & 9.99  \\
2s3\_3n20\_2k & 10,800    & $28.7^{}$   & 9      & 3,356.19 & 4,319.54 & 3,721.06  & 69.86  & N/A   \\
2s3\_3n20\_3k & 10,800    & $12.78^{}$  & 37     & 3,673.57 & 4,143.22 & 643.78    & 48.96  & 13.90 \\
2s3\_3n20\_4k & 10,800    & $12.7^{}$   & 53     & 3,997.06 & 4,504.80 & 308.93    & 48.83  & 15.70 \\
2s3\_3n20\_5k & 10,800    & $8.83^{1}$  & 37     & 4,281.65 & 4,659.52 & 576.60    & 41.97  & 11.83 \\
2s3\_3n25\_2k & 10,800    & N/A         & 1      & N/A      & N/A      & 10,810.55 & N/A    & N/A   \\
2s3\_3n25\_3k & 10,800    & $23.63^{}$  & 11     & 3,624.49 & 4,480.96 & 2,435.22  & 61.12  & N/A   \\
2s3\_3n25\_4k & 10,800    & $18.73^{}$  & 17     & 3,929.73 & 4,665.70 & 1,224.06  & 53.72  & N/A   \\
2s3\_3n25\_5k & 10,800    & $13.25^{}$  & 29     & 4,193.59 & 4,749.30 & 1,194.46  & 47.09  & 13.37 \\
2s3\_4n5\_2k  & 41.19     & $0^{}$      & 135    & 2,127.65 & 2,127.65 & 1.44      & 46.52  & 19.39 \\
2s3\_4n5\_3k  & 4,242.45  & $0^{}$      & 14,933 & 2,709.56 & 2,709.56 & 0.83      & 65.43  & 23.82 \\
2s3\_4n5\_4k  & 2,561.16  & $0^{}$      & 14,803 & 2,857.01 & 2,857.01 & 0.71      & 53.57  & 16.88 \\
2s3\_4n5\_5k  & 65.72        & $0^{}$      & 457    & 3,339.38 & 3,339.38 & 0.43      & 59.64  & 24.28 \\
2s3\_4n10\_2k & 8,216.36  & $0^{}$      & 595    & 3,059.48 & 3,059.48 & 29.40     & 51.15  & 6.04  \\
2s3\_4n10\_3k & 10,800    & $2.77^{}$   & 1,435  & 3,419.82 & 3,514.48 & 27.71     & 61.98  & 12.38 \\
\textbf{2s3\_4n10\_4k} & \textbf{10,146.92} & $\boldsymbol{0^{3}}$     & \textbf{1,775}  & \textbf{3,772.86} & \textbf{3,772.86} & 6.79      & \textbf{58.25}  & \textbf{12.97} \\
2s3\_4n10\_5k & 10,800    & $5.32^{1}$  & 2,727  & 4,061.30 & 4,277.32 & 10.06     & 66.54  & 21.63 \\
2s3\_4n15\_2k & 10,800    & $13^{}$     & 25     & 3,287.61 & 3,714.98 & 1,078.98  & 56.34  & 13.73 \\
2s3\_4n15\_3k & 10,800    & $8.02^{}$   & 37     & 3,779.72 & 4,082.98 & 494.77    & 48.11  & 9.58  \\
2s3\_4n15\_4k & 10,800    & $5.09^{}$   & 51     & 4,103.93 & 4,312.72 & 428.04    & 38.87  & 6.88  \\
2s3\_4n15\_5k & 10,800    & $2.09^{}$   & 93     & 4,456.22 & 4,549.38 & 222.72    & 28.15  & 4.39  \\
2s3\_4n20\_2k & 10,800    & $26.39^{3}$ & 3      & 3,455.62 & 4,367.59 & 6,470.46  & 57.10  & N/A   \\
2s3\_4n20\_3k & 10,800    & $14.38^{}$  & 25     & 3,752.86 & 4,292.36 & 949.23    & 43.19  & 14.97 \\
2s3\_4n20\_4k & 10,800    & $17.63^{}$  & 39     & 4,271.65 & 5,024.80 & 487.14    & 49.47  & 20.69 \\
2s3\_4n20\_5k & 10,800    & $15.63^{3}$ & 45     & 4,735.34 & 5,475.41 & 596.62    & 45.63  & 18.46 \\
2s3\_4n25\_2k & 10,800    & N/A         & 1      & N/A      & N/A      & 10,810.16 & N/A    & N/A   \\
2s3\_4n25\_3k & 10,800    & $39.91^{}$  & 5      & 3,831.72 & 5,361.07 & 4,938.98  & 77.82  & N/A   \\
2s3\_4n25\_4k & 10,800    & $33.26^{}$  & 9      & 4,236.64 & 5,645.85 & 2,789.55  & 69.07  & N/A   \\
2s3\_4n25\_5k & 10,800    & $11.56^{}$  & 23     & 4,621    & 5,154.97 & 1,210.66  & 42.48  & 12.83 \\
2s3\_5n5\_2k  & 93.53     & $0^{}$      & 393    & 1,534    & 1,534    & 1.45      & 71.52  & 7.88  \\
2s3\_5n5\_3k  & 53.49     & $0^{}$      & 289    & 1,629.78 & 1,629.78 & 0.79      & 57.21  & 7.25  \\
2s3\_5n5\_4k  & 13.29     & $0^{}$      & 83     & 1,730.43 & 1,730.43 & 0.63      & 49.69  & 0.45  \\
2s3\_5n5\_5k  & 30.40    & $0^{}$      & 245  & 1921.56 & 1921.56 & 0.22      & 52.74  & 15.60  \\
2s3\_5n10\_2k & 10,800    & $1.73^{}$   & 1,007  & 2,694.37 & 2,740.88 & 43.16     & 51.13  & 10.85 \\
2s3\_5n10\_3k & 10,800    & $6.08^{1}$  & 1,235  & 2,939.82 & 3,118.42 & 29.50     & 57.64  & 14.76 \\
2s3\_5n10\_4k & 10,800    & $9.69^{}$   & 1,935  & 3,169.60 & 3,476.65 & 16.05     & 67.16  & 22.63 \\
2s3\_5n10\_5k & 8,047.76  & $0^{}$      & 2,461  & 3,441.89 & 3,441.89 & 8.23      & 58.58  & 26    \\
2s3\_5n15\_2k & 10,800    & $3.06^{1}$  & 37     & 2,994.03 & 3,085.59 & 1,013.50  & 39.69  & 4.03  \\
2s3\_5n15\_3k & 10,800    & $4.49^{3}$  & 75     & 3,395.46 & 3,547.98 & 219.13    & 45.75  & 7.53  \\
2s3\_5n15\_4k & 10,800    & $10.95^{3}$ & 93     & 3,560.69 & 3,950.44 & 96.04     & 60.38  & 14.57 \\
2s3\_5n15\_5k & 10,800    & $7.55^{}$   & 345    & 4,035.53 & 4,340.15 & 43.19     & 57.78  & 17.21 \\
2s3\_5n20\_2k & 10,800    & $23.2^{}$   & 3      & 3,027.60 & 3,729.99 & 5,325.94  & 57.97  & N/A   \\
2s3\_5n20\_3k & 10,800    & $18.06^{}$  & 27     & 3,472.10 & 4,099.32 & 1,219.99  & 55.50  & 18.77 \\
2s3\_5n20\_4k & 10,800    & $14.45^{}$  & 37     & 3,840.05 & 4,395.09 & 726.75    & 49.70  & 15.90 \\
2s3\_5n20\_5k & 10,800    & $9.4^{}$    & 73     & 4,120.27 & 4,507.46 & 311.77    & 44.26  & 12.99 \\
2s3\_5n25\_2k & 10,800    & $20.7^{}$   & 7      & 3,227    & 3,895    & 4,544.11  & 62.62  & N/A   \\
2s3\_5n25\_3k & 10,800    & $25.67^{}$  & 11     & 3,465.76 & 4,355.51 & 2,503.72  & 71.45  & N/A   \\
2s3\_5n25\_4k & 10,800    & $20.12^{}$  & 21     & 3,838.62 & 4,610.93 & 1,050.47  & 67.86  & 23.90 \\
2s3\_5n25\_5k & 10,800    & $13.76^{}$  & 49     & 4,287.55 & 4,877.31 & 542.75    & 61.73  & 23.12 \\* \bottomrule\\

\rlap{\textsuperscript{1} variation \#1 with parameters: $\kappa = 8$}\\
\rlap{\textsuperscript{2} variation \#2 with parameters: local depth-first search, halfpoint$ = 0.05$}\\
\rlap{\textsuperscript{3} variation \#3 with parameters: local depth-first search, $\kappa = 7$ }\\
\end{longtable}

\end{APPENDICES}


\bibliographystyle{informs2014trsc} 
\bibliography{reference.bib} 


\end{document}